\documentclass[a4paper,11pt,number,preprint]{article}
\usepackage{tikz}
\usetikzlibrary{decorations.pathreplacing}
\usepackage{amsmath,amssymb,amsthm}
\usepackage{cite}
\usepackage{graphicx}
\usepackage{subfigure}
\usepackage{xcolor}
\usepackage{hyperref}
\usepackage{float,placeins}
\usepackage{setspace}
\allowdisplaybreaks[4]
\newtheorem{thm}{Theorem}[section]

 \newtheorem{lem}{Lemma}[section]
 \newtheorem{prop}{Proposition}[section]
 \newtheorem{defn}{Definition}[section]
\newtheorem{rem}{Remark}[section]
\def\N{\mathcal{N}}

\def\J{J_\varepsilon}

\def\d{\partial}
\def\ddj{\dot \Delta_j}

\def\tilde{\widetilde}
\def\hat{\widehat}

\def \zr1{$z_{R,1}$}
\def \zr2{z_{R,2}}
\def \zi1{z_{I,1}}
\def \zi2{z_{I,2}}

\def\W{{\mathcal W}}

\renewcommand{\div}{\mbox{\rm div}\;\!}

\newcommand\R{\mathbb{R}}

\newcommand\Z{\mathbb{Z}}

\def\B{\mathbb{B}}

\theoremstyle{plain}
\begingroup
\theoremstyle{plain}
\newtheorem{theorem}{Theorem}[section]

\newtheorem{lemma}[theorem]{Lemma}
\theoremstyle{definition}

\theoremstyle{remark}

\def\J{{J^{\varepsilon}}}

\def\H0{H_{0}}

\renewcommand{\div}{\mbox{\rm div}\;\!}

\begin{document}
\title{The compressible Euler system with damping in hybrid Besov spaces: global well-posedness and relaxation limit}

\author{Timothée Crin-Barat and Zihao Song\footnote{Corresponding author: szh1995@nuaa.edu.cn}}
\date{}

\maketitle
\begin{abstract}
We investigate the global well-posedness of the compressible Euler system with damping in $\R^d$ ($d \ge 1$) and its relaxation limit toward the porous medium equation. In \cite{CBD3}, the first author and Danchin studied these two problems in hybrid Besov spaces, where the high-frequency components of the solution are bounded in $L^2$-based norms, while the low-frequency components are controlled in $L^p$-based norms with $p \in [2,\max\{4,2d/(d-2)\}]$. Motivated by the observation that the limit system is well-posed in $L^p$-based spaces for $p \in [2,\infty)$, we extend the low-frequency analysis to this full range, thereby providing a more unified framework for studying such relaxation limits.

The core of our proof consists in establishing refined product and commutator estimates describing sharply the interactions between the high, medium, and low-frequency regimes. 
A key observation underlying our analysis is that the product of two functions localized at low frequencies generates only interactions between low and medium frequencies, never purely high-frequency ones. Consequently, for a suitable choice of frequency threshold, the high-frequency projection of the product of two functions localized low frequencies vanishes.

% This novel analysis of nonlinear interactions is widely applicable and can be used to study the whole class of partially dissipative systems satisfying the (SK) condition.
\end{abstract}
\smallbreak

\vspace{2mm}


\section{Introduction}
\subsection{Presentation of the model}
We are concerned with the Cauchy problem of the multi-dimensional compressible Euler system with damping in $\R^d$, for $d\geq1$, which reads
\begin{equation}\label{Euler relax}
\left\{\begin{aligned}
 &\partial_t\rho+\operatorname{div}(\rho u) =0, \\
&\partial_t(\rho u)+\operatorname{div}(\rho u\otimes u)+\nabla P+\frac{1}{\varepsilon}\rho u  =0,\\
&(\rho, u)\big|_{t=0}=(\rho_0, u_0),
\end{aligned} \right.
\end{equation}
where, for $x\in\R^d$ and $t\geq0$, $\rho=\rho(x,t)\in \R$ represents the density, $u=u(x,t)\in\R^d$ the velocity and $\varepsilon>0$ the time-relaxation parameter. We assume that the pressure satisfies the $\gamma$-law:
\begin{eqnarray*}
P(\rho)=A\rho^{\gamma}
\end{eqnarray*}
for some constants $ \gamma>1$ and $A>0.$

In this paper, we focus on strong solutions of system \eqref{Euler relax} associated to initial data that are close to a constant equilibrium $(\bar{\rho},0)$ with $\bar{\rho}>0$. Without loss of generality, we assume that $\bar{\rho}=1$ in the rest of the paper.
We are also concerned with the relaxation limit under the so-called diffusive rescaling:   \begin{equation}\label{DiffusiveRescaling}(\rho^\varepsilon,u^\varepsilon)(t,x)\triangleq (\rho,\varepsilon^{-1}u)(\varepsilon^{-1}t,x).\end{equation}
When expressed in the diffusive scaling \eqref{DiffusiveRescaling}, system \eqref{Euler relax} takes the form
\begin{equation}\label{Eulereps}
\left\{
\begin{array}
[c]{l}%
\partial_{t}\rho^{\varepsilon}+\div(\rho^{\varepsilon} u^{\varepsilon})=0,\\[1ex]
\varepsilon^2 \partial_{t}(\rho^{\varepsilon}u^{\varepsilon})+\varepsilon^2\textrm{div} (\rho^{\varepsilon}u^{\varepsilon}\otimes u^{\varepsilon})+\nabla P(\rho^{\varepsilon})+\rho^{\varepsilon}u^{\varepsilon}=0.
\end{array}
\right.
\end{equation}
As $\varepsilon \to 0$, we expect $(\rho^\varepsilon,u^\varepsilon)$ to converge to $(\mathcal{N},V)$, where $(\mathcal{N},V)$ solves the porous medium system
\begin{equation}\label{PMeq}
\left\{
\begin{array}{l}
\partial_t \mathcal{N}-\Delta P(\mathcal{N})=0,\\[1ex]
V=-\dfrac{\nabla P(\mathcal{N})}{\mathcal{N}},
\end{array}
\right.
\end{equation}
the second equation being Darcy's law.

% and that the above limit system, being parabolic,  
% is well-posed in $L^p$-type spaces.

\subsection{Literature review}

System (\ref{Euler relax}) belongs to the class of partially dissipative hyperbolic systems. Such multi-dimensional first order $n$-component systems read
\begin{equation}\label{hyperbolic}
\partial_t U+\sum^{d}_{k=1}A^k(U)\partial_{x_k} U+\frac{1}{\varepsilon}LU=0,
\end{equation}
where $U=U(x,t)\in \R^n$ is the unknown, the  symmetric matrices valued maps 
$A^k$  ($k=1,\cdots,d$) are assumed to be linear and $L$ is a positive symmetric matrix.

If $L=0$, then (\ref{hyperbolic}) reduces to a classical conservation law. In this case, the global existence theory is open and, as demonstrated in the well-known works of Lax \cite{L} and Serre \cite{Serre}, shocks may develop even for smooth and small initial data. A natural line of investigation then consists in considering hyperbolic system with appropriate source terms, 
providing additional dissipative structures to (\ref{hyperbolic}).

In the pioneering works of Kawashima \cite{K} and of Shizuta and Kawashima \cite{SK}, a sufficient structural condition, known as the (SK) condition, was identified to ensure global well-posedness for small initial data for a class of hyperbolic–parabolic systems. Since then, this theory has been extensively developed. We refer to Yong \cite{Y}, Kawashima and Yong \cite{KY} and Beauchard and Zuazua \cite{BZ} for the proof of existence of global-in-time solutions and the study of their large-time dynamics in Sobolev spaces $H^s$ with $s > \frac{d}{2} + 1$. Later, Kawashima and Xu extended the theory to lower regularity settings and considered initial data in the inhomogeneous Besov space $\dot{\B}^{d/2+1}_{2,1}$. More recently, the first author and Danchin \cite{CBD1,CBD2} analyzed partially dissipative hyperbolic systems in critical homogeneous Besov spaces.
% and their study was further extended to hybrid $L^2-L^p$-based spaces in \cite{CBD3}.

Concerning the compressible Euler system with damping, Chen, Levermore, and Liu \cite{CLL} first introduced appropriate entropy structures to capture the dissipative behavior of the system. The global existence of smooth solutions was then obtained by Yong \cite{Y} under certain restrictions on the entropy. Subsequently, Sideris, Thomases, and Wang \cite{STW}, and Wang and Yang \cite{WY}, established its global well-posedness and justified the asymptotic behavior in Sobolev spaces $H^s$ with $s > \frac{d}{2} + 1$. Xu and Wang \cite{XW} later refined these results by investigating the relaxation limit under critical regularity assumptions.

As for the singular diffusive limit for quasilinear hyperbolic systems, see Marcati and Milani \cite{MM}, Marcati, Milani, and Secchi \cite{MMS}, and Marcati and Rubino \cite{MR} for weak convergence results. In \cite{Junca}, Junca and Rascle were able to justify the relaxation process from the damped Euler equations to the porous media equation in the one-dimensional setting for large global-in-time $BV$ solution and to provide an explicit rate of convergence. Their approach is based on a stream function technique which is related to the mass Lagrangian coordinates. More recently, following the same approach as in \cite{Junca}, Peng et al. in \cite{Peng2019} justified the convergence of partially dissipative hyperbolic systems to parabolic systems globally-in-time in one space dimension and derived a convergence rate of the relaxation process. Using similar techniques, Liang and Shuai in \cite{LiangShuai} generalized the previous result to the multi-dimensional periodic setting (in $\mathbb{T}^3$). Recently, the first author, Peng and Shou \cite{CBPS} recovered a similar result in $\R^3$.\newline
\indent Concerning approaches based on standard energy estimates, Coulombel and Goudon \cite{CG} and Coulombel and Lin \cite{CoulombelLin} justified the relaxation limit in Sobolev spaces, and Xu and Wang \cite{XW} further lowered the regularity requirement on the initial data and justified the limit in inhomogeneous Besov spaces. In \cite{CBD3}, the first author and Danchin justified a strong convergence result providing an explicit rate of convergence for this type of relaxation limit in the multi-dimensional setting. \\ \indent  We also mention further papers related to the investigation of the stability of partially dissipative systems around constant equilibria \cite{Sroczinski25,CBLSZ,BianchiniJinXin,CSZ,JX,CBD1,FayeRodrigues23}.

% Here, we build upon the method developed in \cite{CBD3} and extend their approach.
% \cite{castro2019}

\subsection{Hybrid \texorpdfstring{$L^2-L^p$}{Lp} well-posedness results and aims of the paper}
Since the work of Brenner \cite{B}, it is known that, in general, hyperbolic systems are ill-posed in $L^p$ spaces for $p\ne 2$. Well-posedness results may only hold if $p=2$ or if the matrices $A^j$ commute with each others as they would then be diagonalizable in a same basis and one ends up with decoupled, up to order $0$ terms, transport equations. 

Nevertheless, in the presence of dissipative operators, it turns out that it is possible to construct solutions whose low-frequencies are bounded in $L^p$-based spaces with $p\ne 2$ and whose high frequencies are bounded in $L^2$ ones. Indeed, in \cite{CBD2}, Danchin and the first author established the well-posedness and relaxation limit for the compressible Euler system with damping in hybrid $L^2-L^p$ homogeneous Besov spaces with $p\ne2$. More precisely, the low-frequency part of initial data is assumed to belong to $\dot{\B}^{\frac{d}{p}}_{p,1}$ under the restriction\begin{align}\label{CBDcondp}
p\in\left[2,\min\{4,\frac{2d}{d-2}\}\right]
\end{align}
and the high-frequency part of the initial data is in\footnote{We recall that the embeddings $\dot \B^{\frac d2+1}_{2,1} \hookrightarrow W^{1,\infty}$ and $\dot \B^{\frac{d}{p}}_{p,1} \hookrightarrow L^\infty $
hold. Consequently, under these assumptions on the initial data, the solutions are bounded in
\(L^1_T(W^{1,\infty})\) and \(L^\infty_T(L^\infty)\).} $\dot{\B}^{\frac d2+1}_{2,1}$. We also mention the survey of \cite{DanchinEMS} where a larger class of partially dissipative hyperbolic systems satisfying the (SK) condition are dealt with and the work of Shou, Xu and Zhang \cite{ShouXuZhang} where the regularity assumptions in low frequencies are improved.

In this paper, our aim is to extend the upper bound on $p$ in \eqref{CBDcondp} to $p < \infty$, providing a global well-posedness  and relaxation theory for a broader and more consistent class of initial data. The consistency being related to the following observations:
\begin{itemize}
\item The limit porous media system \eqref{PMeq} is well posed in $L^p$ spaces for all $p < \infty$ (see Proposition~\ref{prop1}).
\item As observed in \cite{CBD3}, the compressible Euler system with damping exhibits a porous-media type behavior in the low-frequency regime.
\end{itemize}
% As a corollary, this also enables us to justify the diffusive relaxation limit in a more consistent framework. 
 % \item The threshold separating the low and high frequencies part of the solution scales as $1/\varepsilon$. Therefore, only the porous media part persists in the limit $\varepsilon \to 0$.

% Our proof consists in developing new product and commutator estimates which describe sharply nonlinear frequency interactions. A key observation that allows us to remove the restriction on $p$ is that when estimating the high-frequency counterpart of a
% product of two functions localized in low frequencies, only terms involving interactions
% between low and medium frequencies arise, never purely high-frequency ones. We refer to Section \ref{sec:strat} for more details concerning these frequency interactions and our strategy of proof.

\section{Main results}
\subsection{Linearization of the system}
Before stating our main results, we present the linearization of the compressible Euler system with damping (\ref{Euler relax}). 
We introduce the rescaled sound speed and velocity: 
\begin{equation}\label{eq:c}  c^\varepsilon(\varepsilon t,x)\triangleq\frac{1}{\check\gamma}\sqrt{\frac{\d P}{\d\rho}}=\frac{(\gamma A)^\frac{1}{2}}{\check\gamma}(\rho(t,x))^{\check\gamma}
\quad \text{and} \quad v^\varepsilon(\varepsilon t,x)=\frac{1}{\varepsilon}u(t,x)\end{equation}
where $ \check{\gamma}=\dfrac{\gamma-1}{2}$.
% Defining the rescaled variables
% \begin{eqnarray}\label{changevariable}
% c^\varepsilon(\varepsilon t,x)=\frac{(4\gamma A)^{\frac{1}{2}}}{\gamma-1}(\rho^{\frac{\gamma-1}{2}}(t,x)-\bar\rho^{\frac{\gamma-1}{2}})\quad \text{and} \quad v^\varepsilon(\varepsilon t,x)=\frac{1}{\varepsilon}u(t,x),
% \end{eqnarray}
% where $c^\varepsilon$ and $v^\varepsilon$ are the diffusively rescaled sound speed and velocity, we obtain
%  \begin{equation}
%  \left\{
%  \begin{array}
%  [c]{l}
% \partial_tc^\varepsilon+v^\varepsilon\cdot\nabla c^\varepsilon+\check\gamma c^\varepsilon\div v^\varepsilon=0,\\[1ex]\displaystyle \varepsilon^2\left(\partial_tv^\varepsilon+v^\varepsilon\cdot \nabla v^\varepsilon\right)+\check\gamma c^\varepsilon\nabla  c^\varepsilon+v^\varepsilon
% =0,
%  \end{array}
%  \right.\label{CEDRelax00}
%  \end{equation}
% where $\check{\gamma}=\frac{\gamma-1}{2}$ and $\bar c=\frac{(4\gamma A)^{\frac{1}{2}}}{\gamma-1}\bar\rho^{\check{\gamma}}$.
Setting  $c:=c^\varepsilon-\bar{c}$ and $v:=v^\varepsilon-0$,
system \eqref{Euler relax} rewrites as
\begin{equation} \left\{ \begin{aligned} &\partial_t c+v\cdot\nabla c+\check{\gamma}(c+\bar{c})\textrm{div}\,v=0,\\ 
&\varepsilon^2(\partial_tv+v\cdot\nabla v)+\check{\gamma}(c+\bar{c})\nabla c+\ v=0, \end{aligned} \right.\label{CED4}
\end{equation} 
where $\bar c=\frac{(4\gamma A)^{\frac{1}{2}}}{\gamma-1}$.
We note that, under this reformulation, system \eqref{CED4} is symmetric. Using the definition of the pressure $P$,  formally, as $\varepsilon\to0$, a solution of \eqref{CED4} is expected to converge to a solution of the porous media equation \eqref{PMeq}.

\subsection{Functional spaces}
% We introduce the Littlewood-Paley decomposition and refer to \cite{HJR} for additional information. Let $\chi$ be a smooth function valued in $[0,1]$ such that $\chi$ is supported in the ball
% $\mathbf{B}(0,\frac{4}{3})=\{\xi\in\mathbb{R}^{d}:|\xi|\leq\frac{4}{3}\}$. We set $\varphi(\xi)=\chi(\xi/2)-\chi(\xi)$ such that $\varphi$
% is supported in the annulus $\mathbf{C}(0,\frac{3}{4},\frac{8}{3})=\{\xi\in\mathbb{R}^{d}:\frac{3}{4}\leq|\xi|\leq\frac{8}{3}\}$ so
% $$\sum_{j\in\mathbb{Z}}\varphi(2^{-j}\xi)=1, \quad \forall\xi\in\mathbb{R}^{d}\backslash\{{0}\}.$$
% For any $j\in \mathbb{Z}$, the homogeneous dyadic blocks $\dot{\Delta}_{j}$ and the low-frequency cut-off operator $\dot{S}_{j}$ are defined by
% $$
% \dot{\Delta}_{j}f:=\mathcal{F}^{-1}(\varphi(2^{-j}\cdot )\mathcal{F}f),\quad\quad \dot{S}_{j}f:= \mathcal{F}^{-1}( \chi (2^{-j}\cdot) \mathcal{F} f),
% $$
% where $\mathcal{F}$ and $\mathcal{F}^{-1}$ stand for the Fourier transform and its inverse. 

% Let $\mathcal{S}_{h}'$ be the set of tempered distributions on $\mathbb{R}^{d}$ such that every $f\in \mathcal{S}_{h}'$ satisfies $f\in \mathcal{S}'$ and $\lim_{j\rightarrow-\infty}\|\dot{S}_{j}f\|_{L^{\infty}}=0$. For $f\in \mathcal{S}_{h}'$, we recall the definition of the homogeneous Besov semi-norms: For any $s\in \R$ and $p\in[1,\infty]$,
% \begin{equation*}
% \|f\|_{\dot{\B}_{p,1}^s}:=\sum_{j\in \Z}2^{js}\|\dot{\Delta}_jf\|_{L^p}.
% \end{equation*}

We begin by recalling the notation associated with the Littlewood–Paley decomposition and Besov spaces. See \cite[Chapter 2]{HJR} for a complete overview. We choose a smooth, radial, non-increasing function $\chi(\xi)$ with compact support in $B(0,\frac{4}{3})$ and $\chi(\xi)=1$ in $B(0,\frac{3}{4})$ such that
$$
\varphi(\xi):=\chi(\frac{\xi}{2})-\chi(\xi),\quad \sum_{j\in \mathbb{Z}}\varphi(2^{-j}\cdot)=1,\quad \text{{\rm{Supp}}}~ \varphi\subset \{\xi\in \mathbb{R}^{d}~|~\frac{3}{4}\leq |\xi|\leq \frac{8}{3}\}.
$$
For any $j\in \mathbb{Z}$, the homogeneous dyadic blocks $\dot{\Delta}_{j}$ and the low-frequency cut-off operator $\dot{S}_{j}$ are defined by
$$
\dot{\Delta}_{j}f:=\mathcal{F}^{-1}(\varphi(2^{-j}\cdot )\mathcal{F}f),\quad\quad \dot{S}_{j}f:= \mathcal{F}^{-1}( \chi (2^{-j}\cdot) \mathcal{F} f),
$$
where $\mathcal{F}$ and $\mathcal{F}^{-1}$ stand for the Fourier transform and its inverse. Throughout the paper, we use the notation $\dot{\Delta}_{j}f:=f_{j}.$

Let $\mathcal{S}_{h}'$ be the set of tempered distributions on $\mathbb{R}^{d}$ such that $\lim\limits_{j\rightarrow-\infty}\|\dot{S}_{j}f\|_{L^{\infty}}=0$. Then, we have
\begin{equation}\nonumber
\begin{aligned}
&f=\sum_{j\in \mathbb{Z}}f_{j}\quad\text{and}\quad \dot{S}_{j}f= \sum_{j'\leq j-1}u_{j'}\quad\text{in}~\mathcal{S}_h'.
\end{aligned}
\end{equation}
the homogeneous Besov space $\dot{\B}^{s}_{p,r}$, for $p,r\in[1,\infty]$ and $s\in \mathbb{R}$, is defined by
$$
\dot{\B}^{s}_{p,r}:=\{f\in \mathcal{S}_{h}'~|~\|f\|_{\dot{\B}^{s}_{p,r}}:=\|\{2^{js}\|uf_{j}\|_{L^p}\}_{j\in\mathbb{Z}}\|_{l^{r}}<\infty\}.
$$

In order to study the system under consideration, we introduce a decomposition of the frequency domain.
We define the frequency threshold 
\begin{align}
    J^\varepsilon=\left\lfloor \textrm{log}_2(\varepsilon^{-1})
\right\rfloor +k_0,
\end{align}
for some $k_0\in \mathbb{Z}$. For any $f\in \mathcal{S}'(\R^d)$, we define 
$$f:=f^{\ell,\varepsilon}+f^{h,\varepsilon}\quad \text{where} \quad  f^{\ell,\varepsilon}:=\sum_{j< \J}\ddj f \quad \text{and} \quad f^{h,\varepsilon}:=\sum_{j\geq \J}\ddj f.$$
Accordingly, we define the frequency-restricted Besov semi-norms
\begin{equation*}
\|f\|_{\dot{\B}_{p,1}^s}^{\ell,\varepsilon}:=\sum_{j<J^\varepsilon}2^{js}\|\dot{\Delta}_jf\|_{L^p} \quad \text{and} \quad
\|f\|_{\dot{\B}_{2,1}^s}^{h,\varepsilon}:=\sum_{j\geq \J}2^{js}\|\dot{\Delta}_jf\|_{L^2}.
\end{equation*}
This decomposition is further refined in Section \ref{sec:moreLP} by introducing additional intermediate frequency regimes.

\subsection{Main results}
We introduce the following functionals
$$X^\varepsilon(t):=X^{\ell,\varepsilon}(t)+X^{h,\varepsilon}(t),$$
where
\begin{eqnarray*}
X^{\ell,\varepsilon}(t):=\| c\|^{\ell}_{L^\infty_T(\dot{\mathbb{B}}_{p,1}^{\frac{d}{p}})\cap L^1_T(\dot{\mathbb{B}}_{p,1}^{\frac{d}{p}+2})}+\varepsilon\| v\|^{\ell}_{L^\infty_T(\dot{\mathbb{B}}_{p,1}^{\frac{d}{p}})}+\| v\|^{\ell}_{L^2_T(\dot{\mathbb{B}}_{p,1}^{\frac{d}{p}})\cap L^1_T(\dot{\mathbb{B}}_{p,1}^{\frac{d}{p}+1})} +\frac{1}{\varepsilon}\|\W\|^{\ell}_{L^1_T(\dot{\B}_{{p},1}^{\frac{d}{{p}}})},
\end{eqnarray*}
with $\W=v+\check\gamma(c+ \bar c)\nabla c$ and
\begin{eqnarray*}
X^{h,\varepsilon}(t):=\varepsilon\| c\|^{h}_{L^\infty_T(\dot{\mathbb{B}}_{2,1}^{\frac{d}{2}+1})}+\frac{1}{\varepsilon}\| c\|^{h}_{L^1_T(\dot{\mathbb{B}}_{2,1}^{\frac{d}{2}+1})}
+\varepsilon^{2}\| v\|^{h}_{L^\infty_T(\dot{\mathbb{B}}_{2,1}^{\frac{d}{2}+1})}+\| v\|^{h}_{L^1_T(\dot{\mathbb{B}}_{2,1}^{\frac{d}{2}+1})}.
\end{eqnarray*}
We also denote $E_0^\varepsilon$ for the functional space of the initial data associated to the norm $X_0^\varepsilon$ defined by
\begin{equation*}\begin{aligned}
X_{0}^\varepsilon:= & \|( c_{0},\varepsilon v_{0})\|^{\ell}_{\dot{\mathbb{B}}_{p,1}^{\frac{d}{p}}}
+\|(\varepsilon  c_{0},\varepsilon^{2} v_{0})\|^{h}_{\dot{\mathbb{B}}_{2,1}^{\frac{d}{2}+1}}.
\end{aligned}\end{equation*}

We first state a global well-posedness result for small initial data.

\begin{thm}\label{thm2}
Let $d\geq1$ and $p\in[2,\infty)$.
There exists a $\eta>0$ and a $k_0\in \mathbb{Z}$ such that for all $\varepsilon>0$, if \begin{align}\label{inithm}\|( c_{0},\varepsilon v_{0})\|^{\ell,\varepsilon}_{\dot{\mathbb{B}}_{p,1}^{\frac{d}{p}}}+\|(\varepsilon  c_{0},\varepsilon^{2} v_{0})\|^{h,\varepsilon}_{\dot{\mathbb{B}}_{2,1}^{\frac{d}{2}+1}}\leq\eta,\end{align}
then the Cauchy problem (\ref{CED4}) admits a unique global-in-time solution $(c,v)$ such that, for all $t>0$, $$X^\varepsilon(t)\leq CX_0^\varepsilon,$$
where $C>0$ is a universal constant.
\end{thm}
\begin{rem}
  In Theorem \ref{thm1}, we remove the restrictions on $d$ and $p$ of \eqref{CBDcondp} appearing in \cite{CBD3}. This allows us to construct global-in-time solutions and to justify the associated relaxation limit for a wider class of initial data.
%   Such a generalization on the pair $(d,p)$
% coincides with the analysis of the linear solution maps, but might be surprising for extending beyond the classical H$\ddot o$lder inequality in nonlinear calculations.
\end{rem}

\begin{rem}
% Although  Theorem \ref{thm2} is stated in a concise way, o 
Our computations in fact yield a more general result, see Theorem \ref{thm1}. In this theorem we define several intermediate frequency regimes in order to precisely describe the frequency interactions and sharpen the assumptions on the initial data. 
  % In particular, the first intermediate regime reflects the observation that, when estimating the high-frequency component of the product of two functions localized at low frequencies, only low and intermediate frequencies appear. 
  These additional frequency regimes are the key ingredients of our proof that enable us to sharpen the upper bound on the space-integrability parameter $p$. See Section \ref{sec:strat} for a detailed description of the proof strategy.
\end{rem}

%   \begin{rem}  In the previous theorem, to simplify the presentation, we defined $J^\varepsilon$ and $\widetilde{J^\varepsilon}$ so that the low-frequency regime encompasses all the medium-frequency regimes that we defined previously. Doing that we lose in precision in the properties satisfied by the solution.
%   In the previous theorem, for the sake of simplicity, we defined $J^\varepsilon$ and $\widetilde{J^\varepsilon}$
%  so that the low-frequency regime includes all the medium-frequency ranges introduced earlier. This choice simplifies the presentation but comes at the cost of losing some precision in the description of the solution’s properties.
% \end{rem}

Next, we show that the solutions constructed in Theorem \ref{thm1}  converges strongly to the solutions of the porous medium equations.

\begin{thm}\label{thm3}
Let $d\geq1$, $p\in[2,\infty)$ and $(c,v)$ be the solution from Theorem \ref{thm2} associated to the initial data $(c_0,v_0)\in E_0$. Let $\mathcal{N}\in\mathcal{C} _b( \mathbb{R} ^+ ; \dot{\mathbb{B}}_{p,1}^{\frac dp})\cap L^1( \mathbb{R}^+;\dot{\mathbb{B}}_{p,1}^{\frac dp+2})$ be the corresponding global solution of \eqref{PMeq}  associated with the initial data $\mathcal{N}_0 \in \dot{\B}^{\frac dp}_{p,1}$ (see Proposition \ref{prop1}). Let $\rho$ and $\rho_0$ be the densities corresponding to $c$ and $c_0$ through the relation \eqref{eq:c} and \eqref{DiffusiveRescaling}. 

Assume that
$$\left\|\rho_{0}^\varepsilon-\mathcal{N}_{0}\right\|_{\mathbb{B}_{p,1}^{\frac{d}{p}-\delta}}\leq C\varepsilon^{\delta}\quad for\,\,\,\delta\in(0,1].$$
Then, as $\varepsilon \rightarrow 0$, we have
$$\rho^\varepsilon-\mathcal{N}\longrightarrow0\quad\textit{strongly in}\quad L^\infty(\mathbb{R}^+;\dot{\mathbb{B}}_{p,1}^{\frac{d}{p}-\delta})\cap L^\frac{2}{1+\delta}(\mathbb{R}^+;\dot{\mathbb{B}}_{p,1}^{\frac{d}{p}+1}),$$
and
$${v}^\varepsilon+\frac{\nabla P({\rho^\varepsilon})}{{\rho}^\varepsilon}\longrightarrow0\quad \text{strongly}\:\,\,in\quad L^1(\mathbb{R}^+;\dot{\mathbb{B}}_{p,1}^{\frac{d}{p}}).$$
Moreover, we have the following quantitative error estimates: for $\delta\in(0,1]$ and $r\in[1,2)$,
$$\|{\rho}^\varepsilon-\mathcal{N}\|_{L^{\infty}(\mathbb{R}^{+};\dot{\mathbb{B}}_{p,1}^{\frac{d}{p}-\delta})\cap
L^{\frac{2}{1+\delta}}(\mathbb{R}^{+};\dot{\mathbb{B}}_{p,1}^{\frac{d}{p}+1})}\leq C\varepsilon^{\delta}\quad \text{and}\quad\left\|{v}^\varepsilon+\frac{\nabla P({\rho}^\varepsilon)}{{\rho}^\varepsilon}\right\|_{L^{r}(\mathbb{R}^{+};\dot{\mathbb{B}}_{p,1}^{\frac{d}{p}})}\leq C\varepsilon.$$
\end{thm}

\begin{rem}
    The rate recovered in Theorem \ref{thm3} depends on the regularity imposed on the initial data. In the case $\delta=1$, we recover the convergence rate from \cite{CBD3} but for a wider range of space-integrability parameter $p$.
\end{rem}

\begin{rem}
The hybrid functional framework and the new frequency-splitting techniques that we develop in this paper may also be extended to general partially dissipative systems of the form~\eqref{hyperbolic} satisfying the (SK) condition. To this end, one would need to combine the arguments presented in the current work with those from~\cite{CBD3,DanchinEMS}.
\end{rem}

\subsection{On hybrid functional frameworks for the compressible Navier-Stokes equations}

\indent Hybrid functional frameworks have been employed to study various types of PDEs, including hyperbolic--parabolic systems such as the compressible Navier--Stokes equations, see Danchin and Charve \cite{NSCLP} and Chen, Miao and Zhang \cite{CMZ}. In their setting, the low-frequency part of the solution is treated in $L^2$-based spaces, while the high-frequency one in $L^p$ ones. An  upper bound similar to condition \eqref{CBDcondp} appears in these works. However, such a condition may not be removed using the techniques developed in the present paper. Indeed, when estimating the low-frequency part of a product of two high-frequency localized functions, terms spanning the entire range of frequencies may arise. Consequently, the upper bound on $p$ seems unavoidable as embeddings of the form $\dot{\B}^{d/p}_{p,1}\hookrightarrow \dot{\B}^{d/2}_{2,1}$ do not hold for $p>2$. 
% In this sense, we believe that the restriction on $p$ for the compressible Navier--Stokes system is sharp. 
Additional remarks on this issue can be found in the Section \ref{sec:strat} and in Remark \ref{remLFHF-NSC}.

We also point out the recent work by Guo, the second author and Yang \cite{GuoSongYang25}, which extends the admissible range of $p$ by exploiting dispersive estimates and assuming stronger regularity assumption on the low-frequency part of the initial data. We believe that an analogous analysis could be performed for the compressible Euler system with damping as dispersive phenomenon are expected in the high-frequency regime.

% \subsubsection{General partially dissipatives systems}

% \subsubsection{Local well-posedness result for the Navier-Stokes equation}
% The classical local well posedness results in \cite{Da}

\subsection{Strategy of the proof}\label{sec:strat}
\subsubsection{Spectral analysis of the linearized system}
We provide a spectral analysis of the linear system associated to \eqref{CED4} for $\check\gamma \bar c=1$:
\begin{equation}\label{lin00}
\left\{
\begin{array}{l}\partial_t{c}+{\div}\,v=0, \\ [1mm]
\varepsilon^2\partial_tv+\nabla{c}+v=0.\\[1mm]
 \end{array} \right.
\end{equation}
Applying the Fourier transform to system \eqref{lin00}, we obtain that $\Omega=(c,v)$ satisfies
\begin{eqnarray*}
\partial_t\hat\Omega+H(\xi)\hat\Omega=0\qquad \text{with} \quad H(\xi)=\begin{pmatrix}
0 & i\xi\\[2ex]
i\xi & \varepsilon^{-1}
\end{pmatrix}.
\end{eqnarray*}
The two eigenvalues of $H(\xi)$ behave asymptotically as
\begin{itemize}
\item For $|\xi|\ll\varepsilon^{-1}$: $\lambda_+\sim\varepsilon^{-1},$ $\lambda_-\sim\varepsilon|\xi|^2$;
\item For $|\xi|\gg\varepsilon^{-1}$: $\Re(\lambda_\pm)\sim\varepsilon^{-1}$.
\end{itemize}
In high frequencies, we observe exponential decay at the rate $1/\varepsilon$ for the whole solution and in low frequencies, one mode is exponentially damped while the other behaves as the heat kernel.

To capture the dissipative behavior of such a partially dissipative system with energy methods in high frequencies, we shall construct Lyapunov functionals in the spirit of the hypocoercivity theory, see \cite{Villani,BZ,NSCL2,CBD1}. Due to the hyperbolic nature of the problem, such energy estimates are restricted to $L^2$-based space to avoid losing a derivative.

In low frequencies, employing a standard hypocoercivity method would not capture the true spectral dynamics, as it would produce a heat-like behavior for both components. A key point is that, for $|\xi|\leq \varepsilon^{-1}$, the eigenvalues are purely real and thus one expects to be able to (partially) decouple the system so as to capture the different behavior expressed by the two modes. Such phenomenon was highlighted in \cite{CBD1,CBD3} and is inspired from the works for Hoff \cite{Hoff} and Haspot \cite{Haspot} related to the introduction of the effective velocity for the compressible Navier-Stokes system. Defining the damped mode $w=v+\nabla c$, system \eqref{lin00} can be rewritten as
	   \begin{equation} \left\{ \begin{aligned} &\partial_tc-\Delta c=-\div w,\\&\varepsilon\partial_tw+\frac{w}{\varepsilon}=-\varepsilon\nabla \Delta c-\varepsilon\nabla \div w. \end{aligned} \right.\label{lind} \end{equation}
      In \eqref{lind}, the left-hand sides captures perfectly the spectral behavior.
As for the right-hand side terms, since they are of higher order, in the low-frequency regime $|\xi|\leq \varepsilon^{-1}$, they can be absorbed by the left-hand side dissipation by means of Bernstein-type inequalities.
In particular, this formulation directly yields $\mathcal{O}(\varepsilon)$ bounds for $w$, which is the key ingredient to justify the relaxation limit.
Moreover, since the dissipative effects are now decoupled, one can estimate both equations in $L^p$-based spaces separately.

\subsubsection{Refined frequency decomposition}
In light of the analysis in the previous section, it is natural to consider an hybrid functional framework where, for some threshold \(J\in\Z\),
\begin{itemize}
    \item \(L^p\) estimates are performed in the low-frequency regime \((|\xi|\leq J)\),
    \item \(L^2\) estimates are performed in the high-frequency regime \((|\xi|\geq J)\).
\end{itemize}
In a linear setting, this approach is valid for any \(p \in [2, \infty)\). In a nonlinear setting, the situation is more complex as interactions between low and high frequencies arise when estimating nonlinearities, requiring additional constraints on the parameter \(p\).  
In what follows, we clarify i) where does such restrictions come from and ii) identify the situations in which they can be relaxed.

\textbf{i) Constraints on $p$.} Let $f,g\in \mathcal{S}_h'(\R^d)$. When one estimates nonlinear terms of the form $fg$ in the low-frequency regime, no additional constraint on $p$ are needed. Indeed, when estimating nonlinear terms such as \(\|f g\|_{\dot{\B}^{d/p}_{p,1}}^\ell\), we use the following decomposition
\begin{align}
    \|f g\|_{\dot{\B}^{d/p}_{p,1}}^\ell\lesssim\|f^\ell g^\ell\|_{\dot{\B}^{d/p}_{p,1}}^\ell+\|f^h g^\ell\|_{\dot{\B}^{d/p}_{p,1}}^\ell+\|f^\ell g^h\|_{\dot{\B}^{d/p}_{p,1}}^\ell+\|f^h g^h\|_{\dot{\B}^{d/p}_{p,1}}^\ell.
\end{align}
Among these terms, the most delicate one is
$\|f^h g^h\|_{\dot{\B}^{d/p}_{p,1}}^\ell$,
since the high-frequency components are controlled in a different functional framework than the low-frequency ones. Nevertheless, this difficulty can be overcome by using embeddings of the form
\begin{align}\label{emb}
\dot{\B}^{d/2}_{2,1} \hookrightarrow \dot{\B}^{d/p}_{p,1},
\qquad \text{for all } p \ge 2,
\end{align}
which allow us to estimate the high-frequency contribution of the nonlinearity in $L^2$-based spaces.
\smallbreak
In the reversed scenario, since the reverse embedding fails, estimating high-frequency terms of the form
\[
\|f^\ell g^\ell\|_{\dot{\B}^{d/2+1}_{2,1}}^h
\]
requires the use of Hölder-type inequalities
\[
L^p\times L^{p^*}\to L^2,
\qquad \frac1p+\frac1{p^*}=\frac12,
\]
combined with the embeddings
\[
\dot{\B}^{\frac{d}{p}-\frac{d}{p^*}}_{p,1}\hookrightarrow L^{p^*},
\qquad p^*=\frac{2p}{p-2},
\qquad p\le p^*.
\]
The condition $p\leq p^*$ yields the constraint $p\leq 4$ in \eqref{CBDcondp}.
\smallbreak
The restriction $p\leq \dfrac{2d}{d-2}$ in \eqref{CBDcondp} stems from regularity considerations in low frequencies.
In that regime, the natural functional framework for our analysis is the space
\(\dot{\B}^{\frac{d}{p}}_{p,1}\).
However, invoking the embedding
\(\dot{\B}^{\frac{d}{p}-\frac{d}{p^*}}_{p,1}\hookrightarrow L^{p^*}\)
leads to a loss of control on the low-frequency components of the solution, since
\[
\frac{d}{p}-\frac{d}{p^*}<\frac{d}{p}.
\]
In practice, this difficulty is addressed by applying the above embedding to the gradient of the
solution rather than to the solution itself.
Consequently, one must control the solution in the space
\(\dot{\B}^{\frac{d}{p}-\frac{d}{p^*}+1}_{p,1}\),
which requires to impose the condition
\[
\frac{d}{p}-\frac{d}{p^*}+1>\frac{d}{p}.
\]
This inequality yields the second constraint in \eqref{CBDcondp}.

\textbf{ii) Removing the constraints.} To relax the restrictions on $p$, we analyze the frequency interactions more precisely. Terms of the form $(f^\ell g^\ell)^h$ can be represented as
\begin{equation*}
( f^\ell g^\ell)^h=\sum_{\substack{2^j\geq J\\2^{j_1},2^{j_2}\leq J}}\ddj(\dot{\Delta}_{j_1}  f\dot{\Delta}_{j_2} g)=\sum_{\substack{2^j\geq J\\2^{j_1},2^{j_2}\leq J}}\int_{\R^d}\varphi\Big(\frac{|\xi|}{2^j}\Big)
\tilde{\varphi}\Big(\frac{|\xi-\eta|}{2^{j_1}})\varphi\Big(\frac{|\eta|}{2^{j_2}}\Big)\hat f(\xi-\eta)\widehat{g}(\eta)d\eta.
\end{equation*}
A key observation is that there exists some fixed constant $1\gg a_0>0$ such that
\begin{eqnarray}\label{obser}\nonumber&&\sum_{\substack{2^j\geq J,\,2^{j_1},\,2^{j_2}\leq J}}\int_{\R^d}\varphi\Big(\frac{|\xi|}{2^j}\Big)
\tilde{\varphi}\Big(\frac{|\xi-\eta|}{2^{j_1}}\Big)\varphi\Big(\frac{|\eta|}{2^{j_2}}\Big)\hat f(\xi-\eta)\widehat{g}(\eta)d\eta\\&=&\sum_{\substack{2^j\geq J,\,2^{j_1},\,2^{j_2}\leq J\\ \max\{2^{j_1},2^{j_2}\}\geq a_0J}}\int_{\R^d}\varphi\Big(\frac{|\xi|}{2^j}\Big)
\tilde{\varphi}\Big(\frac{|\xi-\eta|}{2^{j_1}}\Big)\varphi\Big(\frac{|\eta|}{2^{j_2}}\Big)\hat f(\xi-\eta)\widehat{g}(\eta)d\eta.\end{eqnarray}
Indeed, if $\max\{2^{j_1},2^{j_2}\}<a_0 J$, thanks to the definition of $\varphi$, we find that
$$|\xi|\in[\frac{3}{4}J,\infty),\quad
|\xi-\eta|\in[0,\frac{8}{3}a_0J)\quad \text{and}\quad
|\eta|\in[0,\frac{8}{3}a_0J).$$
Employing the triangle inequality, this implies that
$$|\xi|\leq\big||\xi-\eta|+|\eta|\big|<\frac{16}{3}a_0 J,$$
which, choosing $a_0$ small enough, contradicts $|\xi|\geq\frac{3}{4}J$.
%From \eqref{obser}, one observes that the high-low-low frequency-interactions of the form $(f^\ell g^\ell)^h$, at least one of the source term ($f$ or $g$) shall only be localized within regime $[a_0 J,J]$, which is away from zero and corresponds to a finite sum of the supports of the factors near $J$.  
% Consequently, this implies that the certain interaction $(f^\ell g^\ell)^h$ are supported in narrow frequency bands, their product may spread over a much larger frequency range.

From \eqref{obser}, we infer that for high--low--low frequency interactions of the form $(f^\ell g^\ell)^h$, at least one of the factors $f$ or $g$ must be localized in the intermediate frequency range $[a_0 J,\, J]$. This interval is bounded away from zero and contains only a finite number of dyadic blocks concentrated around the scale $J$. Consequently, the product $(f^\ell g^\ell)^h$ cannot generate arbitrarily high frequencies as it is supported in a relatively narrow frequency band away from zero.

The equality (\ref{obser}) suggests introducing an intermediate medium-frequency regime around the frontier of the low and high-frequency regimes which collects the frequencies generated by nonlinear interactions that escape the high-frequency zone without covering the full low-frequency range.
More concretely, we will consider a medium-frequency regime where the solution is controlled in a \(L^{p_1}\)-based framework with \(p_1\in(2,p)\).  
Nonlinear analysis can then be performed by first applying Hölder inequalities between the medium and high-frequency regimes, yielding a restriction on the intermediate index \(p_1\). Then, performing a similar analysis between the medium and low-frequency regimes will produce a milder restriction on \(p\). This argument can then be iterated by introducing a sufficiently large number of medium-frequency regimes to ultimately reach any \(p\in[2,\infty)\).

In what follows, we first illustrate this mechanism by showing how to reach the case \(p=6\).
For simplicity, we disregard regularity issues and focus solely on the Lebesgue framework and Hölder-type arguments.

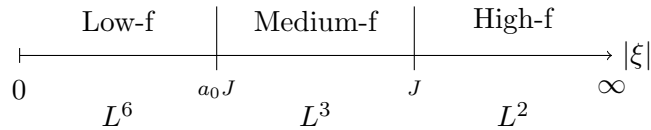
\begin{figure}[h]
\centering
\begin{tikzpicture}[scale=1.3]
\vspace{5mm}
\usetikzlibrary{decorations.pathreplacing}

% Axis
\draw[|->] (0,0) -- (6,0) node[right] {$|\xi|$};

% Division points
\foreach \x in {2,4}{
    \draw (\x,0.5) -- (\x,-0.15);
}

% Labels under axis: 0 and infinity
\node[below] at (0,-0.15) {$0$};
\node[below] at (6,-0.15) {$\infty$};

% Small-font tick labels
\node[below, font=\scriptsize] at (2,-0.15) {$a_0J$};
\node[below, font=\scriptsize] at (4,-0.15) {$J$};

% Top labels for each segment
\node at (1,0.35) {Low-f};
\node at (3,0.35) {Medium-f};
\node at (5,0.35) {High-f};

% Bottom labels for each segment
\node at (1,-0.6) {$L^6$};
\node at (3,-0.6) {$L^{3}$};
\node at (5,-0.6) {$L^2$};

\end{tikzpicture}
\caption{Frequency decomposition diagram for $p=6$.} 
\label{figp6}
\end{figure}

$\bf{Case\: p=6}:$ In this case, the simplest choice is to insert one medium-frequency regime where the solution is bounded in $L^{p_1}$ with $p_1=3$. For some constant $a_0$, the frequency decomposition is represented in Figure \ref{figp6}. 
% The low-frequency regime corresponds to $|\xi|\le a_0 J$, the intermediate-frequency regime to $a_0 J \le |\xi| \le J$, and the high-frequency regime to $|\xi| \ge J$.
\begin{itemize}
    \item In the high-frequency regime, for a fixed $J$, one can show that
$\|f^\ell g^\ell\|_{\dot{\B}^{s}_{2,1}}^h=0$, following the argument used to obtain (\ref{obser}). Hence, the nonlinear term can be decomposed as
$$\|f g\|^h_{\dot{\B}^{s}_{2,1}}=\|f^\ell g^m \|^h_{\dot{\B}^{s}_{2,1}}+\|f^m g^\ell\|^h_{\dot{\B}^{s}_{2,1}}+\mathrm{Remainders}.$$
% \dot{\B}^{h}[f^m,g^\ell]+\dot{\B}^{m}[f^\ell,g^\ell]+\mathrm{Remainders}.$$
For the first two terms, we follow the usual approach used in \cite{CBD1,CBD3}. For instance, we have
$$\|f^\ell g^m \|^h_{\dot{\B}^{s}_{2,1}}\lesssim\|f \|^\ell_{L^{p^*_1}}\| g \|^m_{\dot{\B}^{s}_{p_1,1}},\quad\frac{1}{2}=\frac{1}{p_1}+\frac{1}{p^*_1}$$
and thus $p^*_1\geq6$. 
% The interactions in the remainders can be handled by the classical $L^2\times L^\infty$ argument and Sobolev embedding.
\item In the medium-frequency regime, we have
$$\|f g\|^m_{\dot{\B}^{s}_{p_1,1}}=\|f^\ell g^\ell \|^m_{\dot{\B}^{s}_{p_1,1}}+\mathrm{Remainders}.$$
Applying Hölder inequality, we obtain
 $$\|f^\ell g^\ell \|^m_{\dot{\B}^{s}_{p_1,1}}\lesssim\|f \|^\ell_{L^{\tilde p^*_1}}\| g \|^\ell_{\dot{\B}^{s}_{6,1}},\quad\frac{1}{p_1}=\frac{1}{6}+\frac{1}{\tilde p^*_1}.$$
 
Therefore, this requires the restriction $\tilde p^*_1\geq6$. Notice that choosing $p_1=3$ gives $p^*_1,\tilde p^*_1\geq6$ and thus the previous Hölder estimates are valid. Now, all the low-frequency term appearing are in $L^p$-based space with $p=6$. 
\item In the low-frequency regime, we recall that using embedding of the type \eqref{emb} is enough to control the nonlinear terms.
\end{itemize}

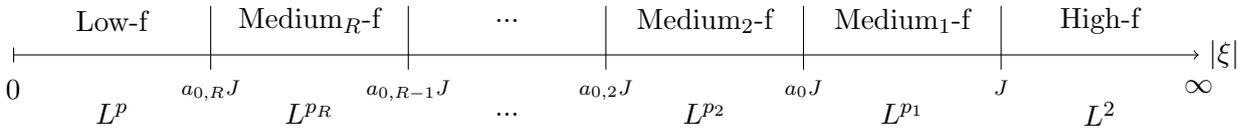
\begin{figure}[h]
\centering
\begin{tikzpicture}[scale=1.3]
\vspace{5mm}
\usetikzlibrary{decorations.pathreplacing}

% Axis
\draw[|->] (0,0) -- (12,0) node[right] {$|\xi|$};

% Division points
\foreach \x in {2,4,6,8,10}{
    \draw (\x,0.5) -- (\x,-0.15);
}

% Labels under axis: 0 and infinity
\node[below] at (0,-0.15) {$0$};
\node[below] at (12,-0.15) {$\infty$};

% Small-font tick labels
\node[below, font=\scriptsize] at (2,-0.15) {$a_{0,R}J$};
\node[below, font=\scriptsize] at (4,-0.15) {$a_{0,R-1}J$};
\node[below, font=\scriptsize] at (6,-0.15) {$a_{0,2}J$};
\node[below, font=\scriptsize] at (8,-0.15) {$a_0J$};
\node[below, font=\scriptsize] at (10,-0.15) {$J$};

% Top labels for each segment
\node at (1,0.35) {Low-f};
\node at (3,0.35) {Medium$_R$-f};
\node at (5,0.35) {...};
\node at (7,0.35) {Medium$_2$-f};
\node at (9,0.35) {Medium$_1$-f};
\node at (11,0.35) {High-f};

% Bottom labels for each segment
\node at (1,-0.6) {$L^{p}$};
\node at (3,-0.6) {$L^{p_R}$};
\node at (5,-0.6) {...};
\node at (7,-0.6) {$L^{p_2}$};
\node at (9,-0.6) {$L^{p_1}$};
\node at (11,-0.6) {$L^{2}$};

\end{tikzpicture}
\caption{Frequency Decomposition Diagram for $p<\infty$.} 
\label{figpinfty}
\end{figure}

$\bf{Case\,p<\infty}:$ If we want to reach larger index $p$, we need to insert additional medium-frequency regimes, each associated with a different space-integrability parameter $(p_i)_{i=1,\ldots,R}$. The first medium-frequency regime collects the frequencies generated by nonlinear interactions that escape the high-frequency zone without covering the full low-frequency range. The second medium-frequency regime collects the frequencies generated by nonlinear interactions that escape the first medium-frequency zone without covering the full low-frequency range. This process is iterated until a sufficient number of medium-frequency regimes have been introduced to achieve the desired integrability. The frequency decomposition is represented in Figure \ref{figpinfty}.

As in the case $p=6$, by choosing $a_{0,1}$ sufficient small, we are able to remove the frequency-interactions of the ``high-low-low" type. Then, choosing some $a_{0,i}$, for $i=2,\ldots,R$, suitably small, we can remove the interactions of the type ``medium-low-low".
We present the analysis of some problematic terms:
\begin{itemize}
\item High-medium$_1$-low interactions:
$$\|f^\ell g^{m_1} \|^h_{\dot{\B}^{s}_{2,1}}\lesssim\|f \|^\ell_{L^{p^*_1}}\| g \|^m_{\dot{\B}^{s}_{p_1,1}},\quad\frac{1}{2}=\frac{1}{p_1}+\frac{1}{p^*_1}.$$

\item Medium$_i$-medium$_{i+1}$-low interactions:
$$\|f^{\ell} g^{m_{i+1}} \|^{m_i}_{\dot{\B}^{s}_{p_i,1}}\lesssim\|f \|^\ell_{L^{p^*_{i+1}}}\| g \|^m_{\dot{\B}^{s}_{p_{i+1},1}},\quad\frac{1}{p_i}=\frac{1}{p_{i+1}}+\frac{1}{p^*_{i+1}}.$$
\item Medium$_R$-low-low interactions:
 $$\|f^\ell g^\ell \|^{m_R}_{\dot{\B}^{s}_{p_R,1}}\lesssim\|f \|^\ell_{L^{ p^*}}\| g \|^\ell_{\dot{\B}^{s}_{p,1}},\quad\frac{1}{p_R}=\frac{1}{p^*}+\frac{1}{p}.$$
\end{itemize} 
In order for the above Hölder inequalities to be valid, we need that the sequence $(p_i)_{i=1,2..R}$ fulfills
\begin{eqnarray}\label{Sequence}
p_1\leq\frac{2p}{p-2},\quad p_{i+1}\leq\frac{p_i p}{p-p_i}\quad \text{and}\quad p_R\leq\frac{p}{2}.
\end{eqnarray}
Under the condition \eqref{Sequence}, we have
\begin{align}
    \|f \|^\ell_{L^{p^*_1}}+\|f \|^\ell_{L^{p^*_{i+1}}}+\|f \|^\ell_{L^{ p^*}} \lesssim \|f \|^\ell_{L^{p}}.
\end{align}
This ensures that all low-frequency terms are controlled in \(L^p\)-based spaces, thereby relaxing the constraint $p\leq 4$ in \eqref{CBDcondp}. As for the condition \(p \leq \frac{2d}{d-2}\) in \eqref{CBDcondp}, the same frequency decomposition, together with a suitable choice of the sequence \((p_i)\) such that
\[
p_{1} \leq \frac{2d}{d - 2}, \quad p_{i+1} \leq \frac{p_i d}{d - p_i} \quad \text{and}\quad p \leq \frac{p_Rd}{d-p_R},
\]
allows us to eliminate this second constraint.

The general class of sequence $(p_i)$ that can be employed is defined in Definition \ref{defp}.

\begin{rem}\label{remLFHF-NSC}
In contrast with the case treated in our work, the analysis for the compressible Navier-Stokes equations is reversed: $L^p$ estimates are performed in the high-frequency regime ($|\xi|\le J$), while the low-frequency regime ($|\xi|\ge J$) is restricted to $L^2$. In this context, the restrictions on $p$ \eqref{CBDcondp}, observed in \cite{CMZ,Haspot}, are due to terms of the form $(f^h g^h)^\ell$, corresponding to
\[
\sum_{\substack{2^j\le J,\, 2^{j_1},2^{j_2}\ge J}} 
\int_{\mathbb{R}^d} 
\varphi\Big(\frac{|\xi|}{2^j}\Big) 
\tilde{\varphi}\Big(\frac{|\xi-\eta|}{2^{j_1}}\Big) 
\varphi\Big(\frac{|\eta|}{2^{j_2}}\Big) 
\hat f(\xi-\eta)\,\widehat{g}(\eta)\,d\eta.
\]
For this term, the Fourier support does not vanish as in our case since one can have $|\xi|\to 0$ and simultaneously $|\xi-\eta|$ and $|\eta|\to\infty$. In other words, when estimating the low-frequency component of a product of two functions localized in high frequencies, contributions spanning the entire frequency range arise. Therefore, we do not expect that using a frequency decomposition analogous to the one presented above to be able to remove the restrictions on $p$ when studying the compressible Navier-Stokes equations.
\end{rem}

\subsection{Outline of the paper}
In Section 3, we state new product and commutator estimates adapted to our hybrid functional framework.
Section 4 and Section 5 are devoted to the proofs of Theorems \ref{thm2} and Theorem \ref{thm3}, respectively.
The proofs of the new product and commutator estimates are collected in the Appendix.

% \begin{rem}
%     \label{remLFHF-NSC}
%     \textcolor{red}{Remark on why our approach fails for Navier-Stokes}
% At the end, we remark that our approach probably could not be applied in extending the $L^2-L^p$ hybrid Besov spaces in terms of the compressible Navier-Stokes equation, which shares the same condition as (\ref{CBDcondp}) (see \cite{D}). In the Navier-Stokes cases, the 
% restriction of $p$ for arises from the corresponding "low-high-high" interactions, which is unavoidable even with aids of ouor medium frequency techinique.
% \end{rem}

\section{Non-classical product and commutator estimates}
In this section, we introduce a new class of product and commutator estimates, which are the cornerstone of our proof.

\subsection{Refined functional spaces}\label{sec:moreLP}
For any $f\in \mathcal{S}'_h(\R^d)$, we define
$$f:=f^{\ell}+\sum\limits^R_{i=1}f^{m_i}+f^{h}$$ with
$$f^{\ell}:=\sum_{j<\J-N_0 R}\ddj f,\quad f^{m_i}:=\sum_{j\in J^\varepsilon_i}\ddj f,\quad f^{h}:=\sum_{j\geq \J}\ddj f,$$
where $J^\varepsilon\in \mathbb{Z}$ corresponds to the threshold between the medium-frequency regime $m_1$ and the high-frequency regime, $R\in \N$ corresponds to the number of medium frequency regimes, $N_0$ is some fixed positive constant and $J^\varepsilon_i$ is defined by
$$J^\varepsilon_i:=\big\{j\in\Z, j\in[\J-N_0 i,\J-N_0 (i-1)[\big\}.$$
We also introduce, for $a,b\in[1,R]$, 
$$f^{[\ell,m_a]}:=\sum_{j<J^\varepsilon-N_0 (a-1)}\ddj z,\quad f^{[m_b,h]}:=\sum_{j\geq \J-N_0 b}\ddj z,\quad f^{[m_a,m_b]}:=\sum_{j\in J^\varepsilon_{[a,b]}}\ddj z,$$
where $J_{[a,b]}$ is defined by
$$J^\varepsilon_{[a,b]}:=\big\{j\in\Z, j\in[\J-N_0 a,\J-N_0(b-1)[\big\}.$$
Accordingly, we define the frequency-restricted Besov semi-norms
\begin{equation*}
\|z\|_{\dot{\B}_{p,1}^s}^{\ell}:=\sum_{j<J^\varepsilon-N_0 R}2^{js}\|\dot{\Delta}_jz\|_{L^p},\quad\|z\|_{\dot{\B}_{{p_i},1}^s}^{m_i}:=\sum_{j\in J^\varepsilon_i}
2^{js}\|\dot{\Delta}_jz\|_{L^{p_i}},
\end{equation*}
\begin{equation*}
\|z\|_{\dot{\B}_{2,1}^s}^{h}:=\sum_{j\geq \J}2^{js}\|\dot{\Delta}_jz\|_{L^2},
\end{equation*}
and
\begin{equation*}
\|z\|_{\dot{\B}_{p,1}^{s}}^{[\ell,m_a]}:=\|z\|_{\dot{\B}_{p,1}^s}^{\ell}+\sum^R_{i=a}\|z\|_{\dot{\B}_{{p_i},1}^s}^{m_i},\quad
\|z\|_{\dot{\B}_{p,1}^{s}}^{[m_a,m_b]}:=\sum^a_{i=b}\|z\|_{\dot{\B}_{p_i,1}^s}^{m_i},
\end{equation*}
\begin{equation*}
\|z\|_{\dot{\B}_{p,1}^s}^{[m_b,h]}:=\sum^b_{i=1}\|z\|_{\dot{\B}_{p_i,1}^s}^{m_i}+\|z\|_{\dot{\B}_{2,1}^s}^{h}.
\end{equation*}
\subsection{Product and commutator laws}
We introduce the product and commutator estimates adapted to our hybrid functional framework.
\begin{lemma}\label{product}
Let $d\geq1$, $s>0$ and $p\in[2,\infty)$. Assume $\{p_i\}_{i\in[1,R]}\in(2,p)$ to be an increasing sequence and let $f,g\in \mathcal{S}'(\R^d)$.
\begin{enumerate}
    \item \textit{High-frequency product law:} For fixed $s,s_1,s_2>0$, if $p_1$ fulfills the following condition
\begin{align}\label{Condition1}
s+\frac{2d}{p_1}-\frac{d}{2}>0,\quad p_1>\frac{2p}{p-2},\quad s-\frac{d}{2}=s_1+s_2-\frac{2d}{p_1},
\end{align}
then there exists an integer $N_0>0$  such that
$$\displaylines{
\sum_{j\geq J^\varepsilon}2^{js}\|\ddj( fg)\|_{L^2}\leq C\Bigl(\|f\|_{\dot{\B}_{p,1}^{\frac dp-\frac {d}{p_1^*}}}^{[\ell,m_1]}\|g\|_{\dot{\B}_{p_1,1}^{s}}^{m_1}
+\|f\|_{L^\infty}\|g\|_{\dot{\B}_{2,1}^{s}}^{h}
+ \|g\|_{\dot{\B}_{p,1}^{\frac dp-\frac {d}{p_1^*}}}^{[\ell,m_1]}\|f\|_{\dot{\B}_{p_1,1}^{s}}^{m_1}\hfill\cr\hfill
+\|g\|_{L^\infty}\|f\|_{\dot{\B}_{2,1}^{s}}^{h}+\|f\|_{\dot{\B}_{p_1,1}^{s_{1}}}^{[m_1,h]}\|g\|_{\dot{\B}_{p_1,1}^{s_{2}}}^{[m_1,h]}\Bigl),}$$
where $\frac{1}{2}=\frac{1}{p_1}+\frac{1}{p^{*}_1}$ and the constant $C$ depends on $ R, s, p, p_i$ and $d$.

\item \textit{Medium-frequency product law:} For fixed $s,s_1,s_2>0$ and $p_i\in (2,p)$, if $p_{i+1}$ fulfills the following condition
\begin{align}\label{Condition2}
s+\frac{2d}{p_{i+1}}-\frac{d}{p_{i}}>0,\quad p_{i+1}>\frac{p_i p}{p-p_i},\quad s-\frac{d}{p_i}=s_1+s_2-\frac{2d}{p_{i+1}}
\end{align}
then there exists an integer $N_0$ such that
$$\displaylines{
\sum_{j\in J_i}2^{js}\|\ddj( fg)\|_{L^{p_i}}\leq C\Bigl(\|f\|_{\dot{\B}_{p,1}^{\frac dp-\frac {d}{p_{i+1}^*}}}^{[\ell,m_{i+1}]}\|g\|_{\dot{\B}_{p_{i+1},1}^{s}}^{m_{i+1}}
+\|f\|_{L^\infty}\|g\|_{\dot{\B}_{p_{i},1}^{s}}^{[m_{i},h]}+\|g\|_{\dot{\B}_{p,1}^{\frac dp-\frac {d}{p_{i+1}^*}}}^{[\ell,m_{i+1}]}\|f\|_{\dot{\B}_{p_{i+1},1}^{s}}^{m_{i+1}}
\hfill\cr\hfill+\|g\|_{L^\infty}\|f\|_{\dot{\B}_{p_{i},1}^{s}}^{[m_{i},h]}+ \|f\|_{\dot{\B}_{p_{i+1},1}^{s_{1}}}^{[m_{i+1},h]}\|g\|_{\dot{\B}_{p_{i+1},1}^{s_{2}}}^{[m_{i+1},h]}\Bigl)}$$
where $\frac{1}{p_i}=\frac{1}{p_{i+1}}+\frac{1}{p^{*}_{i+1}}$.
\end{enumerate}
\end{lemma}

\begin{lemma}\label{commutator}
Let $d\geq1$, $s>0$ and $p\in[2,\infty)$. Assume $(p_i)_{i\in[1,R]}\in(2,p)$ to be an increasing sequence and let $f,g\in \mathcal{S}'(\R^d)$. For $j\in\Z,$ denote $\mathfrak{R}_j\triangleq \dot S_{j-1}f\,\ddj g-\ddj(fg)$.

\begin{enumerate}
    \item High-frequency commutator estimates: For fixed $s,s_1,s_2>0$, if $p_1$ fulfills \eqref{Condition1}, then there exists a constant $C$ depending only on $s,$ $p$ and $d$ such that
$$\displaylines{
\sum_{j\geq J^\varepsilon}2^{js}\|\mathfrak{R}_j\|_{L^2}\leq C\Bigl(\|\nabla f\|_{\dot{\B}_{p,1}^{\frac dp-\frac {d}{p_1^*}}}^{[\ell,m_1]}\|g\|_{\dot{\B}_{p_1,1}^{s-1}}^{m_1}
+\|\nabla f\|_{L^\infty}\|g\|_{\dot{\B}_{2,1}^{s-1}}^{h}+
\|g\|_{\dot{\B}_{p,1}^{\frac dp-\frac {d}{p_1^*}}}^{[\ell,m_1]}\|f\|_{\dot{\B}_{p_1,1}^{s}}^{m_1}\hfill\cr\hfill
+\|g\|_{L^\infty}\|f\|_{\dot{\B}_{2,1}^{s}}^{h}+ \|g\|_{\dot{\B}_{p_1,1}^{s_{1}}}^{[m_1,h]}\|f\|_{\dot{\B}_{p_1,1}^{s_{2}}}^{[m_1,h]}\Bigl).}$$

\item  Medium-frequency commutator estimates: For fixed $s,s_1,s_2>0$ and $p_i\in (2,p)$, if $p_{i+1}$ fulfills the following (\ref{Condition2}), then there exists  a constant $C$ depending only on $s,$ $p$ and $d$ such that
$$\displaylines{
\sum_{j\in J_i}2^{js}\|\mathfrak{R}_j\|_{L^{p_i}}\leq C\Bigl(\|\nabla f\|_{\dot{\B}_{p,1}^{\frac dp-\frac {d}{p_{i+1}^*}}}^{[\ell,m_{i+1}]}\|g\|_{\dot{\B}_{p_{i+1},1}^{s-1}}^{m_{i+1}}
+\|\nabla f\|_{L^\infty}\|g\|_{\dot{\B}_{p_{i},1}^{s-1}}^{[m_{i},h]}+
\|g\|_{\dot{\B}_{p,1}^{\frac dp-\frac {d}{p_{i+1}^*}}}^{[\ell,m_{i+1}]}\|f\|_{\dot{\B}_{p_{i+1},1}^{s}}^{m_{i+1}}\hfill\cr\hfill
+\|g\|_{L^\infty}\|f\|_{\dot{\B}_{p_{i},1}^{s}}^{[m_{i},h]}+ \|g\|_{\dot{\B}_{p_{i+1},1}^{s_1}}^{[m_{i+1},h]}\|f\|_{\dot{\B}_{p_{i+1},1}^{s_2}}^{[m_{i+1},h]}\Bigl).}$$
\end{enumerate}
\end{lemma}
The proof of the above lemmas are given in Appendix \ref{Appendix}.

\section{Proof of Theorem \ref{thm2}}\label{sec:gwp}

This section is devoted to the proof of Theorem \ref{thm2}. In fact, our analysis yields the following more general statement which implies directly Theorem \ref{thm2}.

\begin{thm}\label{thm1}
Let $d\geq1$, $p\in[2,\infty)$, $\varepsilon>0$. Assume that $(c_0,v_0)\in E_0$.
% and assume that for every $$( c^{\ell,\J}_{0},\varepsilon v^{\ell,\J}_{0})\in{\dot{\mathbb{B}}_{p,1}^{\frac{d}{p}}}, \quad( c^{m_{i},\J}_{0},\varepsilon v^{m_{i},\J}_{0})\in{\dot{\mathbb{B}}_{p_i,1}^{\frac{d}{p_i}}} \quad \text{and} \quad( \varepsilon c^{h,\J}_{0},\varepsilon^{2}v^{h,\J}_{0})\in{\dot{\mathbb{B}}_{2,1}^{\frac{d}{2}+1}}.$$ 
There exists a $\eta>0$, a $k_0\in\mathbb{Z}$, a $R>0$ and an increasing sequence $(p_i)_{i=1,\ldots,R}$ (depending on $d$ and $p$) such that if \begin{align}\label{inithm2}X_0\leq\eta,\end{align} where $X_0$ is given in (\ref{initial}), then the Cauchy problem (\ref{CED4}) admits a unique global-in-time solution $(c,v)$ satisfying
\begin{eqnarray*}
&&  (c,\varepsilon v)^{\ell,\J}\in{L^\infty_T(\dot{\mathbb{B}}_{p,1}^{\frac{d}{p}})}, \,\,(\nabla c,v)^{\ell,\J}\in{ L^1(\mathbb{R}^+;\dot{\mathbb{B}}_{p,1}^{\frac{d}{p}+1})},\,\,\frac{1}{\varepsilon}\W^{\ell}\in{L^1_T(\dot{\B}_{{p},1}^{\frac{d}{{p}}})} \\
 &&  (c,\varepsilon v)^{m_i,\J}\in{L^\infty_T(\dot{\mathbb{B}}_{p_i,1}^{\frac{d}{p_i}})},\,\,(\nabla c,v)^{m_i,\J}\in{ L^1_T(\dot{\mathbb{B}}_{p_i,1}^{\frac{d}{p_i}+1})},\,\,\frac{1}{\varepsilon}\W^{m_i}\in{L^1_T(\dot{\B}_{{p},1}^{\frac{d}{{p}}})}\\
 &&  (\varepsilon c,\varepsilon^2 v)^{h,\J}\in{L^\infty_T(\dot{\mathbb{B}}_{2,1}^{\frac{d}{2}+1})} \quad \text{and} \quad (\frac{1}{\varepsilon} c, v)^{h,\J}\in{ L^1_T(\dot{\mathbb{B}}_{2,1}^{\frac{d}{2}+1})}.
\end{eqnarray*}
Furthermore, we have the uniform estimate
$$X(t)\leq CX_0,$$
where $X(t)$ is defined in \eqref{Xfull} and $C>0$ is a universal constant.
\end{thm}
In Theorem \ref{thm1}, we introduced the following functional:
\begin{align}\label{Xfull}
X(t):=X^\ell(t)+\sum^R_{i=1}X^{m_i}(t)+X^h(t),
\end{align}
where
\begin{eqnarray*}
X^\ell(t):=\| c\|^{\ell}_{L^\infty_T(\dot{\mathbb{B}}_{p,1}^{\frac{d}{p}})\cap L^1_T(\dot{\mathbb{B}}_{p,1}^{\frac{d}{p}+2})}+\varepsilon\| v\|^{\ell}_{L^\infty_T(\dot{\mathbb{B}}_{p,1}^{\frac{d}{p}})}+\| v\|^{\ell}_{L^2_T(\dot{\mathbb{B}}_{p,1}^{\frac{d}{p}})\cap L^1_T(\dot{\mathbb{B}}_{p,1}^{\frac{d}{p}+1})} +\frac{1}{\varepsilon}\|\W\|^{\ell}_{L^1_T(\dot{\B}_{{p},1}^{\frac{d}{{p}}})},
\end{eqnarray*}
\begin{eqnarray*}
X^{m_i}(t):=\| c\|^{m_i}_{L^\infty_T(\dot{\mathbb{B}}_{p_i,1}^{\frac{d}{p_i}})\cap L^1_T(\dot{\mathbb{B}}_{p_i,1}^{\frac{d}{p_i}+2})}+\varepsilon\| v\|^{m_i}_{L^\infty_T(\dot{\mathbb{B}}_{p_i,1}^{\frac{d}{p_i}})}+\| v\|^{m_i}_{L^2_T(\dot{\mathbb{B}}_{p_i,1}^{\frac{d}{p_i}})\cap L^1_T(\dot{\mathbb{B}}_{p_i,1}^{\frac{d}{p_i}+1})}+\frac{1}{\varepsilon}\|\W\|^{m_i}_{L^1_T(\dot{\B}_{{p_i},1}^{\frac{d}{{p_i}}})}
\end{eqnarray*}
and
\begin{eqnarray*}
X^{h}(t):=\varepsilon\| c\|^{h}_{L^\infty_T(\dot{\mathbb{B}}_{2,1}^{\frac{d}{2}+1})}+\frac{1}{\varepsilon}\| c\|^{h}_{L^1_T(\dot{\mathbb{B}}_{2,1}^{\frac{d}{2}+1})}
+\varepsilon^{2}\| v\|^{h}_{L^\infty_T(\dot{\mathbb{B}}_{2,1}^{\frac{d}{2}+1})}+\| v\|^{h}_{L^1_T(\dot{\mathbb{B}}_{2,1}^{\frac{d}{2}+1})}.
\end{eqnarray*}
We also denotes $E_0$ for the functional space associated to the norm $X_0$ defined by
\begin{equation}\label{initial}\begin{aligned}
X_{0}:= & \|( c_{0},\varepsilon v_{0})\|^{\ell}_{\dot{\mathbb{B}}_{p,1}^{\frac{d}{p}}}+\sum^R_{i=1}\|( c_{0},\varepsilon v_{0})\|^{m_i}_{\dot{\mathbb{B}}_{p_i,1}^{\frac{d}{p_i}}}
+\|(\varepsilon  c_{0},\varepsilon^{2} v_{0})\|^{h}_{\dot{\mathbb{B}}_{2,1}^{\frac{d}{2}+1}}.
\end{aligned}\end{equation}

\begin{rem}
  Compared with Theorem \ref{thm2}, Theorem \ref{thm1} imposes weaker integrability assumptions on the medium-frequency components of the initial data. In Theorem \ref{thm2}, the entire medium-frequency regime is assumed to belong to $L^2$, whereas in Theorem \ref{thm1} these components are only required to lie in $L^{p_i}$ spaces, for a suitable sequence $(p_i)$.
\end{rem}

\begin{rem} 
The class of admissible sequences $(p_i)_{i=1,\ldots,R}$ that we employ in Theorem \ref{thm1} is specified in Definition~\ref{defp}. 
Its construction plays a central role in extending the admissible range of the space–integrability index $p$. 
Each element $p_i$ corresponds to the integrability index associated with one of the medium–frequency regimes inserted between the low and high-frequency ones. 
These intermediate scales allows us to progressively improve the integrability range by means of Hölder-type estimates. 

% We stress that the condition~\eqref{inithm} is weaker than the one used in~\cite{CBD3}, where it was assumed that $c^{m_i,\J}\in \dot{\mathbb{B}}_{p,1}^{\frac{d}{p}}$ for all $i=1,\ldots,R$. 
% This can be seen from the embedding $\dot{\mathbb{B}}_{q_1,1}^{\frac{d}{q_1}} \hookrightarrow \dot{\mathbb{B}}_{q_2,1}^{\frac{d}{q_2}}$ for $q_1 \le q_2$.
\end{rem}
\begin{rem} 
In Theorem~\ref{thm1}, the condition on the initial data~\eqref{inithm} depends on the choice of the sequence~$(p_i)$. It can be simplified by imposing additional regularity in the medium-frequency regime. Theorem \ref{thm2} illustrates this simplification and is a direct consequence of Theorem~\ref{thm1}.
\end{rem}

The main component of the proof of Theorem \ref{thm1} relies on  a uniform priori estimate, which will be presented in the next subsection.

\subsection{A priori estimates}
We begin by introducing the notion of admissible sequences. This concept corresponds to the choice of medium-frequency space-integrability parameters that can be selected to extend the range of the low-frequency integrability parameter $p$.
\begin{defn}\label{defp}
Let $p\in(2,\infty)$, $d\geq 1$ and $R\in\Z$ such that $R>1$. We define the set of admissible sequence $\mathcal{S}_{p,R}$ as follows:\\
If $p\in (2,4]$,
\begin{eqnarray*}
\!\!\!\!&&\mathcal{S}_{p,R}=\Big\{(p_i)_{i=1,\ldots,R}\in(2,p), p_i<p_{i+1};\:
p_1\leq\frac{2d}{d-2}, \, p_{i+1}\leq\frac{p_i d}{d-p_i},\, p\leq\frac{p_{R}d}{d-p_{R}}\Big\}. \quad\qquad
\end{eqnarray*}
If $p\in (4,\infty)$,
\begin{eqnarray*}
\!\!\!\!&&\mathcal{S}_{p,R}=\Big\{(p_i)_{i=1,\ldots,R}\in(2,p), p_i<p_{i+1};\: p_1\leq\min\{\frac{2p}{p-2},\frac{2d}{d-2}\},\\&&\nonumber
 \begin{array}
 [c]{l}
\qquad\qquad\qquad p_{i+1}\leq\min\{\frac{p_i p}{p-p_i},\frac{p_id}{d-p_i}\},\,\,p\leq\min\{2p_R,\frac{p_{R}d}{d-p_{R}}\}\end{array}\Big\}.
\end{eqnarray*}
\end{defn}
Here is an example of sequence satisfying \eqref{Sequence}:
$$(p_i)_{i\in[1,R]}=\{p_i=2+\frac{4i}{p-2},\,\, R=\left[\frac{(p-2)(p-4)}{8}\right]+1\}.$$

The main goal of this subsection is to prove the following priori estimates:
\begin{lem}\label{priori}
Let $d\geq1$ and $p\in(2,\infty)$. Suppose that $(p_i)\in\mathcal{S}_{p,R}$ and that $(c,v)$ is a smooth solution of system (\ref{CED4}). The following a priori estimate holds true
\begin{eqnarray}\label{estpriori}
X(t)\lesssim X_0+X^2(t)+X^3(t).
\end{eqnarray}
\end{lem}

The rest of this subsection is devoted to the proof of Lemma \ref{priori}.

\subsubsection{High frequencies estimates}
Let $j\geq J^\varepsilon$ and assume that $\check{\gamma}=1$. Applying the localization operator $\ddj$ to the system \eqref{CED4}, we obtain
\begin{equation*}
\left\{
\begin{array}{l}\partial_t{c}_{j}+\dot S_{j-1}v\cdot\nabla{c}_{j}+{\div}v_{j}=R_j^1+R_j^2-\dot S_{j-1}c\div{v}_{j}, \\ [1mm]
\varepsilon^2\partial_tv_{j}+\varepsilon^2\dot S_{j-1}v\cdot\nabla v_{j}+\nabla{c}_{j}+v_{j}=R_j^3+R_j^4-\dot S_{j-1}c\nabla c_{j},\\[1mm]
 \end{array} \right.
\end{equation*}
where
\begin{eqnarray*}R_j^1:=\dot S_{j-1}v\cdot\nabla\ddj c-\ddj(v\cdot\nabla c),\quad R_j^2:=\dot S_{j-1}c\div\ddj v-\ddj(c\div v),\end{eqnarray*}
\begin{eqnarray*}R_j^3:=\varepsilon^2\big(\dot S_{j-1}v\cdot\nabla\ddj v-\ddj(v\cdot\nabla v)\big)\quad \text{and} \quad R_j^4:=\dot S_{j-1}c\nabla\ddj c-\ddj(c\nabla c).\end{eqnarray*}
Taking the $L^2$-inner product with 
$(c_j,v_j)$, we obtain
\begin{multline}\label{E1}
\frac{1}{2}\frac{d}{dt}\left\|({c}_{j},\varepsilon v_{j})\right\|_{L^2}^2+\frac{1}{\varepsilon^2}\left\|\varepsilon v_{j}\right\|_{L^2}^2\lesssim\left\|\nabla v\right\|_{L^\infty}\left\|({c}_{j},\varepsilon v_{j})\right\|_{L^2}^2\\
+\left\|(R_j^1,R_j^2)\right\|_{L^2}\left\|{c}_{j}\right\|_{L^2}+\left\|\frac{1}{\varepsilon}( R_j^3,R_j^4)\right\|_{L^2}\left\|\varepsilon v_{j}\right\|_{L^2},
\end{multline}
where we used that
\begin{eqnarray*}
\int\dot S_{j-1}v\cdot\nabla{c}_{j}\, {c}_{j}dx\lesssim\left\|\nabla v\right\|_{L^\infty}\left\|{c}_{j}\right\|_{L^2}^2
\end{eqnarray*}
and
\begin{eqnarray*}
\varepsilon^{2}\int\dot S_{j-1}v\cdot\nabla{v}_{j}\,{v}_{j}dx\lesssim\left\|\nabla v\right\|_{L^\infty}\left\|\varepsilon{v}_{j}\right\|_{L^2}^2.
\end{eqnarray*}
To recover dissipation for the density, we compute the time-derivative of the following cross terms:
\begin{multline}\label{E2}
2^{-2j}\frac{d}{dt}\langle\nabla c_j,v_j\rangle+\frac{1}{\varepsilon^{2}}\left\| c_{j}\right\|_{L^2}^2\lesssim\left\| v_{j}\right\|_{L^2}^2+2^{-2j}\frac{1}{\varepsilon^{2}}\langle\nabla c_j,v_j\rangle\\+2^{-j}\frac{1}{\varepsilon}\left\|Q^1_{j}\right\|_{L^2}\left\|\varepsilon v_{j}\right\|_{L^2}+2^{-j}\frac{1}{\varepsilon^{2}}\left\|Q^2_{j}\right\|_{L^2}\left\|{c}_{j}\right\|_{L^2},
\end{multline}
where $Q^{1}=v\cdot \nabla c+c\cdot\nabla v$ and  $Q^2=\varepsilon^2v\cdot \nabla v-c\cdot\nabla c$.
Using the fact $2^j\gtrsim 1/\varepsilon$ implies that $1\lesssim\varepsilon 2^j$, we have
$$\left\| v_{j}\right\|_{L^2}^2=\frac{1}{\varepsilon^2}\left\|\varepsilon v_{j}\right\|_{L^2}^2,\quad
2^{-2j}\frac{1}{\varepsilon^{2}}\langle\nabla c_j,v_j\rangle\lesssim \frac{1}{\varepsilon}\|c_j\|_{L^2}\|v_j\|_{L^2}\lesssim\frac{1}{\varepsilon^{2}}(\|c_j\|^{2}_{L^2}+\|\varepsilon v_j\|^{2}_{L^2}).$$
Hence, multiplying (\ref{E2}) by a constant $\eta>0$ sufficiently small and adding it to (\ref{E1}), we obtain
\begin{multline}\label{E3}
\frac{1}{2}\frac{d}{dt}\mathcal{L}_j+\frac{1}{\varepsilon^2}\left\|({c}_{j},\varepsilon v_{j})\right\|_{L^2}^2\lesssim\left\|\nabla v\right\|_{L^\infty}\left\|({c}_{j},\varepsilon v_{j})\right\|_{L^2}^2
+\left\|(R_j^1,R_j^2)\right\|_{L^2}\left\|{c}_{j}\right\|_{L^2}\\+\left\|\frac{1}{\varepsilon}( R_j^3,R_j^4)\right\|_{L^2}\left\|\varepsilon v_{j}\right\|_{L^2}+2^{-j}\frac{1}{\varepsilon}\left\|Q^1_{j}\right\|_{L^2}\left\|\varepsilon v_{j}\right\|_{L^2}+2^{-j}\frac{1}{\varepsilon^{2}}\left\|Q^2_{j}\right\|_{L^2}\left\|{c}_{j}\right\|_{L^2},
\end{multline}
where \begin{align}
    \mathcal{L}_j=\left\|({c}_{j},\varepsilon v_{j})\right\|_{L^2}^2+2^{-2j}\langle\nabla c_j,v_j\rangle \sim \left\|({c}_{j},\varepsilon v_{j})\right\|_{L^2}^2
\end{align}
since $2^{-2j}\langle\nabla c_j,v_j\rangle\lesssim\|c_j\|^{2}_{L^2}+\|\varepsilon v_j\|^{2}_{L^2}.$

Using Lemma \ref{SimpliCarre}, multiplying by $2^{(\frac{d}{2}+1)j}$ and summing on $j\geq \J$, we obtain
\begin{equation*}\begin{aligned}
\|({c},\varepsilon v)\|_{\dot{\B}_{2,1}^{\frac{d}{2}+1}}^h  +\frac{1}{\varepsilon^2}\|({c},\varepsilon v)\|_{L^1_T(\dot{\B}_{2,1}^{\frac{d}{2}+1})}^h \leq
X_{0}+C\int_0^t\|\nabla v\|_{L^\infty}\|({c},\varepsilon v)\|_{\dot{\B}_{p,1}^{\frac{d}{2}+1}}^h\\
   +\int_0^t\big(\sum_{j\geq \J}2^{(\frac{d}{2}+1)j}\left\|(R_j^1,R_j^2,\frac{1}{\varepsilon} R_j^3,\frac{1}{\varepsilon}R_j^4)\right\|_{L^2}+\sum_{j\geq \J}2^{\frac{d}{2}j}\left\|(\frac{1}{\varepsilon} Q^1_{j},\frac{1}{\varepsilon^{2}}Q^2_{j})\right\|_{L^2}\big),
\end{aligned}
\end{equation*}
which, multiplied by $\varepsilon$, gives
\begin{multline*}
X^h(t)\lesssim X_{0}+\varepsilon\int_0^t\|\nabla v\|_{L^\infty}(\|c\|_{\dot{\B}_{p,1}^{\frac{d}{2}+1}}^h+\varepsilon\|v\|_{\dot{\B}_{p,1}^{\frac{d}{2}+1}}^h)\\
+\varepsilon\int_0^t\big(\sum_{j\geq \J}2^{(\frac{d}{2}+1)j}\left\|(R_j^1,R_j^2,\frac{1}{\varepsilon} R_j^3,\frac{1}{\varepsilon}R_j^4)\right\|_{L^2}+\sum_{j\geq \J}2^{\frac{d}{2}j}\left\|(\frac{1}{\varepsilon} Q^1_{j},\frac{1}{\varepsilon^{2}}Q^2_{j})\right\|_{L^2}\big).
\end{multline*}
Concerning the nonlinearities, we first claim the following inequalities:
\begin{eqnarray}\label{v infty}
\|\nabla v\|_{L^1_T(L^\infty)}\lesssim X(t),\quad \varepsilon^{2}\|\nabla v\|_{L^\infty_T( L^\infty)}\lesssim X(t), \quad\varepsilon\|\nabla a\|_{L^\infty_T( L^\infty)}\lesssim X(t).
\end{eqnarray}
We provide the estimate only for the first inequality, as the other terms can be handled with analogous calculations. We have
\begin{eqnarray}
\|\nabla v\|_{L^\infty}\lesssim\|v\|^{\ell}_{\dot{\mathbb{B}}_{p,1}^{\frac{d}{p}+1}}+\sum^R_{i=1}\|v\|^{m_i}_{\dot{\mathbb{B}}_{p_i,1}^{\frac{d}{p_i}+1}}
+\|v\|^{h}_{\dot{\mathbb{B}}_{2,1}^{\frac{d}{2}+1}},
\end{eqnarray}
which, integrated-in-time, leads to the desired inequality. We thus obtain
\begin{multline}\label{pp1}
\varepsilon\int_0^t\|\nabla v\|_{L^\infty}(\|c\|_{\dot{\B}_{p,1}^{\frac{d}{2}+1}}^h+\varepsilon\|v\|_{\dot{\B}_{p,1}^{\frac{d}{2}+1}}^h)
\leq\Big(\|v\|_{L^1_T(\dot{\B}_{p,1}^{\frac{d}{p}+1})}^{\ell}
+\sum^R_{i=1}\|v\|_{L^1_T(\dot{\B}_{p_{i},1}^{\frac{d}{p_{i}+1}})}^{m_i}+\|v\|_{L^1_T(\dot{\B}_{2,1}^{\frac{d}{2}+1})}^{h}\Big)
\\ \cdot\big(\varepsilon\|c\|_{L^\infty_T(\dot{\B}_{2,1}^{\frac{d}{2}+1})}^h
+\varepsilon^{2}\|v\|_{L^\infty_T(\dot{\B}_{2,1}^{\frac{d}{2}+1})}^h\big)\lesssim X^2(t).
\end{multline}
Hence we are left with commutators $R^i_j$ and products $Q^i_j$.

{\bf Estimate for $R_j^1$ and $R_j^2$}. In fact, taking $s=\frac{d}{2}+1$ in Lemma \ref{commutator} yields
\begin{multline*}
\sum_{j\geq \J}2^{(\frac{d}{2}+1)j}\|R_j^1\|_{L^2}\leq C\Bigl(\|\nabla v\|_{\dot{\B}_{p,1}^{\frac dp-\frac {d}{p_1^*}}}^{[\ell,m_1]}\|\nabla c\|_{\dot{\B}_{p_1,1}^{\frac{d}{2}}}^{m_1}
+\|\nabla v\|_{L^\infty}\|\nabla c\|_{\dot{\B}_{2,1}^{\frac{d}{2}}}^{h}
\hfill\cr\hfill+
\|\nabla c\|_{\dot{\B}_{p,1}^{\frac dp-\frac {d}{p_1^*}}}^{[\ell,m_1]}\|v\|_{\dot{\B}_{p_1,1}^{\frac{d}{2}+1}}^{m_1}
+\|\nabla c\|_{L^\infty}\|v\|_{\dot{\B}_{2,1}^{\frac{d}{2}+1}}^{h}+ \|\nabla c\|_{\dot{\B}_{p_1,1}^{\frac{d}{p_1}}}^{[m_1,h]}\|v\|_{\dot{\B}_{p_1,1}^{\frac{d}{p_1}+1}}^{[m_1,h]}\Bigl).
\end{multline*}
selecting $p_1$ such that for $r=1-\frac{d}{2}+\frac{d}{p_1}$ fulfills $r>0$, there naturally holds
\begin{eqnarray*}
\|\nabla (c,v)\|_{\dot{\B}_{p,1}^{\frac dp-\frac {d}{p_1^*}}}^{[\ell,m_1]}\lesssim\varepsilon^{r}\big(\|(c,v)\|_{\dot{\B}_{p,1}^{\frac dp}}^{\ell}+
\sum^R_{i=1}\| (c,v)\|_{\dot{\B}_{p_i,1}^{\frac {d}{p_i}}}^{m_i}\big),
\end{eqnarray*}
and
\begin{eqnarray*}
\|\nabla c\|_{\dot{\B}_{p_1,1}^{\frac{d}{2}}}^{m_1}\lesssim\varepsilon^{-r}\| c\|_{\dot{\B}_{p_1,1}^{\frac{d}{p_1}+2}}^{m_1},\quad
\| v\|_{\dot{\B}_{p_1,1}^{\frac{d}{2}+1}}^{m_1}\lesssim\varepsilon^{1-r}\| v\|_{\dot{\B}_{p_1,1}^{\frac{d}{p_1}+1}}^{m_1},
\end{eqnarray*}
which implies
\begin{multline}
\varepsilon\int_0^t\|\nabla v\|_{\dot{\B}_{p,1}^{\frac dp-\frac {d}{p_1^*}}}^{[\ell,m_1]}
\|\nabla c\|_{\dot{\B}_{p_1,1}^{\frac{d}{2}}}^{m_1}\lesssim\varepsilon\Big(\| v\|_{L^\infty_T(\dot{\B}_{p,1}^{\frac dp})}^{\ell}+
\sum^R_{i=1}\| v\|_{L^\infty_T(\dot{\B}_{p_i,1}^{\frac {d}{p_i}})}^{m_i}\Big)
\| c\|_{L^1_T(\dot{\B}_{p_1,1}^{\frac{d}{p_1}+2})}^{m_1} \lesssim X^2(t).
\end{multline}
The second term can be handled in the same way as we did to obtain (\ref{pp1}). For the third term, we get
\begin{multline}
\varepsilon\int_0^t\|\nabla c\|_{\dot{\B}_{p,1}^{\frac dp-\frac {d}{p_1^*}}}^{[\ell,m_1]}
\|v\|_{\dot{\B}_{p_1,1}^{\frac{d}{2}+1}}^{m_1}\lesssim\varepsilon^{2}\Big(\|c\|_{L^\infty_T(\dot{\B}_{p,1}^{\frac dp})}^{\ell}+
\sum^R_{i=1}\|c\|_{L^\infty_T(\dot{\B}_{p_i,1}^{\frac {d}{p_i}})}^{m_i}\Big)
\|v\|_{L^1_T(\dot{\B}_{p_1,1}^{\frac{d}{p_1}+1})}^{m_1} \lesssim \varepsilon X^2(t).
\end{multline}
For the fourth term, utilizing (\ref{v infty}), it holds
\begin{eqnarray}
\varepsilon\int_0^t\|\nabla c\|_{L^\infty}\|v\|_{\dot{\B}_{2,1}^{\frac{d}{2}+1}}^h\leq\varepsilon\|\nabla c\|_{L^\infty_T( L^\infty)}\|v\|_{L^1_T(\dot{\B}_{2,1}^{\frac{d}{2}+1})}^h \lesssim X^2(t).
\end{eqnarray}
Finally, for the last term, we get
\begin{multline*}
\varepsilon\int_0^t\|v\|_{\dot{\B}_{p_1,1}^{\frac{d}{p_1}+1}}^{[m_1,h]}\|\nabla c\|_{\dot{\B}_{p_1,1}^{\frac{d}{p_1}}}^{[m_1,h]}\lesssim\varepsilon\Big(\| v\|_{L^2_T(\dot{\B}_{p_1,1}^{\frac{d}{p_1}+1})}^{m_1}\| c\|_{L^2_T(\dot{\B}_{p_1,1}^{\frac{d}{p_1}+1})}^{m_1}+\| v\|_{L^2_T(\dot{\B}_{p_1,1}^{\frac{d}{p_1}+1})}^{m_1}\| c\|_{L^2_T(\dot{\B}_{2,1}^{\frac{d}{2}+1})}^{h}\\
\| v\|_{L^2_T(\dot{\B}_{2,1}^{\frac{d}{2}+1})}^{h}\| c\|_{L^2_T(\dot{\B}_{p_1,1}^{\frac{d}{p_1}+1})}^{m_1}+
\| v\|_{L^2_T(\dot{\B}_{2,1}^{\frac{d}{2}+1})}^{h}\| c\|_{L^2_T(\dot{\B}_{2,1}^{\frac{d}{2}+1})}^{h}\Big)\lesssim X^2(t),
\end{multline*}
where we used
$$\varepsilon\| v\|_{L^2_T(\dot{\B}_{p_1,1}^{\frac{d}{p_1}+1})}^{m_1}\lesssim\| v\|_{L^2_T(\dot{\B}_{p_1,1}^{\frac{d}{p_1}})}^{m_1}
\lesssim X(t).$$
Therefore we get
\begin{equation}\label{R1}
\begin{aligned}
\varepsilon\int^t_0\sum_{j\geq \J}2^{(\frac{d}{2}+1)j}\|R_j^1\|_{L^2}\lesssim X^2(t).
\end{aligned}
\end{equation}
The commutator $R_j^2$ can be handled similarly.

{\bf Estimates for $R_j^3$}. One has
\begin{multline*}
\sum_{j\geq \J}2^{(\frac{d}{2}+1)j}\|R_j^3\|_{L^2}\leq C\Bigl(\|\nabla v\|_{\dot{\B}_{p,1}^{\frac dp-\frac {d}{p_1^*}}}^{[\ell,m_1]}\|\nabla v\|_{\dot{\B}_{p_1,1}^{\frac{d}{2}}}^{m_1}
+\|\nabla v\|_{L^\infty}\|\nabla v\|_{\dot{\B}_{2,1}^{\frac{d}{2}}}^{h}+ \|\nabla v\|_{\dot{\B}_{p_1,1}^{\frac{d}{p_1}}}^{[m_1,h]}\|v\|_{\dot{\B}_{p_1,1}^{\frac{d}{p_1}+1}}^{[m_1,h]}\Bigl).
\end{multline*}
Using that
\begin{multline*}
\varepsilon^{2}\int_0^t\|\nabla v\|_{\dot{\B}_{p,1}^{\frac dp-\frac {d}{p_1^*}}}^{[\ell,m_1]}
\|\nabla v\|_{\dot{\B}_{p_1,1}^{\frac{d}{2}}}^{m_1}\lesssim\varepsilon\Big(\| v\|_{L^\infty_T(\dot{\B}_{p,1}^{\frac dp})}^{\ell}+
\sum^R_{i=1}\| v\|_{L^\infty_T(\dot{\B}_{p_i,1}^{\frac {d}{p_i}})}^{m_i}\Big)
\|v\|_{L^1_T(\dot{\B}_{p_1,1}^{\frac{d}{p_1}+1})}^{m_1}\lesssim X^2(t),
\end{multline*}
\begin{eqnarray*}
\varepsilon^{2}\int_0^t\|\nabla v\|_{L^\infty}\|v\|_{\dot{\B}_{2,1}^{\frac{d}{2}+1}}^h\leq\varepsilon^{2}\|\nabla v\|_{L^\infty_T( L^\infty)}\|v\|_{L^1_T(\dot{\B}_{2,1}^{\frac{d}{2}+1})}^h \lesssim X^2(t),
\end{eqnarray*}
\begin{eqnarray*}
\varepsilon^{2}\int_0^t\|v\|_{\dot{\B}_{p_1,1}^{\frac{d}{p_1}+1}}^{[m_1,h]}\|\nabla v\|_{\dot{\B}_{p_1,1}^{\frac{d}{p_1}}}^{[m_1,h]}\lesssim\varepsilon^{2}\big(\| v\|_{L^2_T(\dot{\B}_{p_1,1}^{\frac{d}{p_1}+1})}^{m_1}+
\| v\|_{L^2_T(\dot{\B}_{2,1}^{\frac{d}{2}+1})}^{h}\big)^{2} \lesssim X^2(t),
\end{eqnarray*}
we arrive at
$$\sum_{j\geq \J}2^{(\frac{d}{2}+1)j}\|R_j^3\|_{L^2}\lesssim X^2(t).$$

{\bf Estimate for $R_j^4$}. We have
\begin{multline*}
\sum_{j\geq J}2^{(\frac{d}{2}+1)j}\|R_j^4\|_{L^2}\leq C\Bigl(\|\nabla a\|_{\dot{\B}_{p,1}^{\frac dp-\frac {d}{p_1^*}}}^{[\ell,m_1]}\|\nabla a\|_{\dot{\B}_{p_1,1}^{\frac{d}{2}}}^{m_1}
+\|\nabla a\|_{L^\infty}\|\nabla a\|_{\dot{\B}_{2,1}^{\frac{d}{2}}}^{h}+ \|\nabla a\|_{\dot{\B}_{p_1,1}^{\frac{d}{p_1}}}^{[m_1,h]}\|a\|_{\dot{\B}_{p_1,1}^{\frac{d}{p_1}+1}}^{[m_1,h]}\Bigl).
\end{multline*}
Using that
\begin{multline*}
\int_0^t\|\nabla a\|_{\dot{\B}_{p,1}^{\frac dp-\frac {d}{p_1^*}}}^{[\ell,m_1]}
\|\nabla a\|_{\dot{\B}_{p_1,1}^{\frac{d}{2}}}^{m_1}\lesssim\Big(\|a\|_{L^\infty_T(\dot{\B}_{p,1}^{\frac dp})}^{\ell}+
\sum^R_{i=1}\| a\|_{L^\infty_T(\dot{\B}_{p_i,1}^{\frac {d}{p_i}})}^{m_i}\Big)
\|a\|_{L^1_T(\dot{\B}_{p_1,1}^{\frac{d}{p_1}+2})}^{m_1} \lesssim X^2(t),
\end{multline*}
\begin{eqnarray*}
\int_0^t\|\nabla a\|_{L^\infty}\|a\|_{\dot{\B}_{2,1}^{\frac{d}{2}+1}}^h\leq\|\nabla a\|_{L^\infty_T( L^\infty)}\|a\|_{L^1_T(\dot{\B}_{2,1}^{\frac{d}{2}+1})}^h \lesssim X^2(t)
\end{eqnarray*}
and
\begin{eqnarray*}
\int_0^t\|a\|_{\dot{\B}_{p_1,1}^{\frac{d}{p_1}+1}}^{[m_1,h]}\|\nabla a\|_{\dot{\B}_{p_1,1}^{\frac{d}{p_1}}}^{[m_1,h]}\lesssim\big(\| a\|_{L^2_T(\dot{\B}_{p_1,1}^{\frac{d}{p_1}+1})}^{m_1}+
\|a\|_{L^2_T(\dot{\B}_{2,1}^{\frac{d}{2}+1})}^{h}\big)^{2} \lesssim X^2(t),
\end{eqnarray*}
we conclude that $R^4_j$ satisfies $$\sum_{j\geq \J}2^{(\frac{d}{2}+1)j}\|R_j^4\|_{L^2}\lesssim X^2(t).$$

{\bf Estimates for $Q^1_{j}$ and $Q^2_{j}$}. We write
\begin{eqnarray}\label{Q1}
Q^1_{j}=\dot S_{j-1}v\cdot\nabla\ddj c+\dot S_{j-1}c\cdot\nabla\ddj v+R_j^1+R_j^2,
\end{eqnarray}
\begin{eqnarray}\label{Q2}
Q^2_{j}=\varepsilon^{2}\dot S_{j-1}v\cdot\nabla\ddj v+\dot S_{j-1}c\cdot\nabla\ddj c+R_j^3+R_j^4.
\end{eqnarray}
We have
\begin{eqnarray*}
\sum_{j\geq \J}2^{\frac{d}{2}j}\|\dot S_{j-1}v\cdot\nabla\ddj c\|_{L^2}\leq C\|v\|_{L^\infty}\| c\|_{\dot{\B}_{2,1}^{\frac{d}{2}+1}}^{h}
\end{eqnarray*}
and
\begin{eqnarray*}
\sum_{j\geq J}2^{\frac{d}{2}j}\|\dot S_{j-1}c\cdot\nabla\ddj v\|_{L^2}\leq C\|c\|_{L^\infty}\| v\|_{\dot{\B}_{2,1}^{\frac{d}{2}+1}}^{h},
\end{eqnarray*}
which imply
\begin{align*}
\int_0^t\| v\|_{L^\infty}\|c\|_{\dot{\B}_{2,1}^{\frac{d}{2}+1}}^h&\leq\frac{1}{\varepsilon}\Big(\varepsilon\|v\|_{L^\infty_T(\dot{\B}_{p,1}^{\frac{d}{p}})}^{\ell}
+\varepsilon\sum^R_{i=1}\|v\|_{L^\infty_T(\dot{\B}_{p_{i},1}^{\frac{d}{p_{i}}})}^{m_i}
+\varepsilon^{2}\|v\|_{L^\infty_T(\dot{\B}_{2,1}^{\frac{d}{2}+1})}^{h}\Big)\|c\|_{L^1_T(\dot{\B}_{2,1}^{\frac{d}{2}+1})}^h\\ &\lesssim X^2(t)
\end{align*}
and
\begin{align*}
\int_0^t\| c\|_{L^\infty}\|v\|_{\dot{\B}_{2,1}^{\frac{d}{2}+1}}^h&\leq\Big(\|c\|_{L^\infty_T(\dot{\B}_{p,1}^{\frac{d}{p}+1})}^{\ell}
+\sum^R_{i=1}\|c\|_{L^\infty_T(\dot{\B}_{p_{i},1}^{\frac{d}{p_{i}}+1})}^{m_i}
+\varepsilon\|c\|_{L^\infty_T(\dot{\B}_{2,1}^{\frac{d}{2}+1})}^{h}\Big)\|v\|_{L^1_T(\dot{\B}_{2,1}^{\frac{d}{2}+1})}^h\\ &\lesssim X^2(t).
\end{align*}
As for remainders, they can be handled in the same fashion as before by noticing that
$$\int_0^t\sum_{j\geq \J}2^{\frac{d}{2}j}\left\|( R^1_{j},R^2_{j})\right\|_{L^2}\lesssim\varepsilon\int_0^t\sum_{j\geq \J}2^{(\frac{d}{2}+1)j}\left\|( R^1_{j},R^2_{j})\right\|_{L^2}\lesssim X^2(t).$$
The terms $Q^1$ and $Q^2$ can be bounded in a similar way. This concludes the high-frequency a priori estimates as we showed that
\begin{eqnarray}\label{priori1}
X^{h}(t) \lesssim X_0+X^2(t).
\end{eqnarray}

\subsubsection{Medium frequencies estimates}
Let $J^\varepsilon-RN_0\leq j< J^\varepsilon$ and $(p_i)\in S_{p,R}$.
Let us recall
\begin{equation} \left\{ \begin{aligned} &\partial_t c+v\cdot\nabla c+\check{\gamma}(c+\bar{c})\textrm{div}\,v=0,\\ 
&\varepsilon^2(\partial_tv+v\cdot\nabla v)+\check{\gamma}(c+\bar{c})\nabla c+\ v=0. \end{aligned} \right.\label{CED40}
\end{equation} 
Defining the damped mode $w= \varepsilon(v+\check{\gamma}(c+\bar{c})\nabla c)$, we have
$$\partial_{t}w+ \frac{1}{\varepsilon^2}w=-(\check{\gamma}\bar c)^2\varepsilon\nabla\div v+Q^3,$$
where
$$Q^3=-\varepsilon v\nabla v-\check{\gamma}\bar c\varepsilon \nabla(v\cdot\nabla c)-(\check{\gamma}\bar c)^2\varepsilon \nabla(c\div v)+\check{\gamma}\varepsilon\partial_t(c\nabla c).$$
Hence, one has
\begin{equation}\label{ww}\begin{aligned}
\|w\|_{\dot{\B}_{{p_i},1}^{\frac{d}{{p_i}}}}^{m_{i} } +\frac{1}{\varepsilon^{2}}\int_0^t\|w\|_{\dot{\B}_{{p_i},1}^{\frac{d}{{p_i}}}}^{m_{i}} \leq
X_{0}+\varepsilon\int_0^t\sum_{j\in J^\varepsilon_i}2^{(\frac{d}{{p_i}}+2)j}\| v_j\|_{L^{p_i}}
+\int_0^t\sum_{j\in J^\varepsilon_i}2^{\frac{d}{{p_i}}j}\|Q^3_{j}\|_{L^{p_i}}.
\end{aligned}
\end{equation}
Using the definition of the damped mode $w$, the density fulfills the following heat equation:
$$\partial_t c-\check{\gamma}\bar c\Delta c=-\frac{1}{\varepsilon}\textrm{div}w+Q^4$$
where
$$Q^4=-v\cdot \nabla c-\check{\gamma}c \div v+\check{\gamma}\div(c\cdot\nabla c).$$
Classical estimates for the heat equation leads to
\begin{equation}\label{cc}\begin{aligned}
\|c\|_{\dot{\B}_{{p_i},1}^{\frac{d}{{p_i}}}}^{m_{i} } +\int_0^t\|c\|_{\dot{\B}_{p_i,1}^{\frac{d}{{p_i}}+2}}^{m_{i}} \leq
X_{0}+\frac{1}{\varepsilon}\int_0^t\|w\|_{\dot{\B}_{p_i,1}^{\frac{d}{{p_i}}}}^{m_{i}}
+\int_0^t\sum_{j\in J^\varepsilon_i}2^{\frac{d}{{p_i}}j}\|Q^4_{j}\|_{L^{p_i}}.
\end{aligned}
\end{equation}
Consequently adding (\ref{ww}) to (\ref{cc}) multiplied by a constant $K$  sufficiently small and choosing $k_0$ sufficiently small, we obtain
\begin{multline}
\|w\|_{L^\infty_T(\dot{\B}_{{p_i},1}^{\frac{d}{{p_i}}})}^{m_{i} } +\frac{1}{\varepsilon^2}\|w\|_{L^1_T(\dot{\B}_{{p_i},1}^{\frac{d}{{p_i}}})}^{m_{i} } +K\big(\|c\|_{L^\infty_T(\dot{\B}_{{p_i},1}^{\frac{d}{{p_i}}})}^{m_{i} }+\|c\|_{L^1_T(\dot{\B}_{{p_i},1}^{\frac{d}{{p_i}}+2})}^{m_{i} }\big)\\ \leq
X_{0}+\| v_j\|_{L^1_T(\dot{\B}_{{p_i},1}^{\frac{d}{{p_i}}})}^{m_{i}}+\int_0^t\sum_{j\in J^\varepsilon_i}2^{\frac{d}{{p_i}}j}\| (Q^3_{j},Q^4_{j})\|_{L^{p_i}}.
\end{multline}
For the component $v$, we recover
\begin{equation}\begin{aligned}
\varepsilon\|v\|_{L^\infty_T(\dot{\B}_{{p_i},1}^{\frac{d}{{p_i}}})}^{m_{i} } \lesssim\|w\|_{L^\infty_T(\dot{\B}_{{p_i},1}^{\frac{d}{{p_i}}})}^{m_{i} }+\varepsilon\|c\|_{L^\infty_T(\dot{\B}_{{p_i},1}^{\frac{d}{{p_i}}})}^{m_{i} },
\end{aligned}
\end{equation}
\begin{equation}\begin{aligned}
\|v\|_{L^2_T(\dot{\B}_{{p_i},1}^{\frac{d}{{p_i}}})}^{m_{i} } \lesssim\frac{1}{\varepsilon}\|w\|_{L^2_T(\dot{\B}_{{p_i},1}^{\frac{d}{{p_i}}})}^{m_{i} }+\|c\|_{L^2_T(\dot{\B}_{{p_i},1}^{\frac{d}{{p_i}}+1})}^{m_{i} }
\end{aligned}
\end{equation}
and
\begin{equation}\begin{aligned}
\|v\|_{L^1_T(\dot{\B}_{{p_i},1}^{\frac{d}{{p_i}}+1})}^{m_{i} } \lesssim\frac{1}{\varepsilon}\|w\|_{L^1_T(\dot{\B}_{{p_i},1}^{\frac{d}{{p_i}}})}^{m_{i} }+\|\nabla c\|_{L^1_T(\dot{\B}_{{p_i},1}^{\frac{d}{{p_i}}+1})}^{m_{i} }.
\end{aligned}
\end{equation}
Consequently using that $2^j\leq\frac{2^{k_0}}{\varepsilon}$ implies that $1\leq 2^{k_0}2^{-j}\varepsilon^{-1}$, choosing $k_0$ appropriately small, we arrive at
\begin{eqnarray}\label{errt}
X^{m_i}(t) \leq
X_{0}+\int_0^t\sum_{j\in J^\varepsilon_i}2^{\frac{d}{{p_i}}j}\| (Q^3_{j},Q^4_{j})\|_{L^{p_i}}.
\end{eqnarray}
Next, we focus on the nonlinear terms.

{\bf Estimate for $Q^3_{j}$}. We first recall the definition of $R^1, R^2, R^3$ and that
\begin{eqnarray}\label{Q3}
\ddj(v\cdot\nabla c)=\dot S_{j-1}v\cdot\nabla\ddj c+R_j^1,\quad \ddj(c\div v)=\dot S_{j-1}c\div\ddj v+R_j^2;
\end{eqnarray}
\begin{eqnarray}
\ddj(v\cdot\nabla v)=\dot S_{j-1}v\cdot\nabla\ddj v+R_j^3.
\end{eqnarray}

Inspired by (\ref{Q1}), we have
\begin{eqnarray}\label{SSDDFF}
\nonumber\int_0^t\sum_{j\in J^\varepsilon_i}2^{\frac{d}{{p_i}}j}\|\varepsilon \nabla(\dot{S}_{j-1}v\cdot\nabla \ddj c)\|_{L^{p_i}}&\lesssim&\int_0^t\sum_{j\in J^\varepsilon_i}2^{\frac{d}{{p_i}}j}\|\dot{S}_{j-1}v\cdot\nabla \ddj c\|_{L^{p_i}}\\
\nonumber&\lesssim&\int_0^t\sup_{j\in J^\varepsilon_i}\|\dot{S}_{j-1}v\|_{L^\infty}\|c\|_{\dot{\B}_{p_i,1}^{\frac{d}{p_i}+1}}^{m_i}\\
&\leq&\|v\|_{L^2 L^\infty}\|c\|_{L^2_T(\dot{\B}_{p_i,1}^{\frac{d}{p_i}+1})}^{m_i}\lesssim X^2(t),
\end{eqnarray}
\begin{eqnarray*}
\int_0^t\sum_{j\in J^\varepsilon_i}2^{\frac{d}{{p_i}}j}\|\varepsilon \nabla(\dot{S}_{j-1}c\div \ddj v)\|_{L^{p_i}}&\lesssim&\int_0^t\sum_{j\in J^\varepsilon_i}2^{\frac{d}{{p_i}}j}\|\dot{S}_{j-1}c\div \ddj v\|_{L^{p_i}}\\&\lesssim&\int_0^t\sup_{j\in J^\varepsilon_i}\|\dot{S}_{j-1}c\|_{L^\infty}\|v\|_{\dot{\B}_{p_i,1}^{\frac{d}{p_i}+1}}^{m_i}\\
&\leq&\|c\|_{L^\infty_T( L^\infty)}\|v\|_{L^1_T(\dot{\B}_{p_i,1}^{\frac{d}{p_i}+1})}^{m_i}\lesssim X^2(t),
\end{eqnarray*}
\begin{eqnarray*}
\varepsilon\int_0^t\sum_{j\in J_i}2^{\frac{d}{{p_i}}j}\|\dot{S}_{j-1}v\cdot\nabla \ddj v\|_{L^{p_i}}&\lesssim&\varepsilon\int_0^t\sup_{j\in J^\varepsilon_i}\|\dot{S}_{j-1}v\|_{L^\infty}\|v\|_{\dot{\B}_{p_i,1}^{\frac{d}{p_i}+1}}^{m_i}\\
&\leq&\varepsilon\|v\|_{L^\infty L^\infty}\|v\|_{L^1_T(\dot{\B}_{p_i,1}^{\frac{d}{p_i}+1})}^{m_i}\lesssim X^2(t).
\end{eqnarray*}

Concerning $R^1_{j}$, using Lemma \ref{commutator}, one obtains
\begin{multline}\label{SSDDFFGG}
\sum_{j\in J_i}2^{\frac{d}{p_i}j}\|\varepsilon \nabla R_j^1\|_{L^{p_i}}\leq\sum_{j\in J_i}2^{\frac{d}{p_i}j}\| R_j^1\|_{L^{p_i}}\leq C\Bigl(\|\nabla v\|_{\dot{\B}_{p,1}^{\frac dp-\frac {d}{p_{i+1}^*}}}^{[\ell,m_{i+1}]}\|\nabla c\|_{\dot{\B}_{p_{i+1},1}^{\frac{d}{p_{i}}-1}}^{m_{i+1}}
+\|\nabla v\|_{L^\infty}\|\nabla c\|_{\dot{\B}_{p_i,1}^{\frac{d}{p_i}-1}}^{[m_{i},h]}\hfill\cr\hfill
+
\|\nabla c\|_{\dot{\B}_{p,1}^{\frac dp-\frac {d}{p_{i+1}^*}}}^{[\ell,m_{i+1}]}\|v\|_{\dot{\B}_{p_{i+1},1}^{\frac{d}{p_i}}}^{m_{i+1}}
+\|\nabla c\|_{L^\infty}\|v\|_{\dot{\B}_{p_i,1}^{\frac{d}{p_i}}}^{[m_{i},h]}+ \|\nabla c\|_{\dot{\B}_{p_{i+1},1}^{\frac{d}{p_{i+1}}}}^{[m_{i+1},h]}\|v\|_{\dot{\B}_{p_{i+1},1}^{\frac{d}{p_{i+1}}}}^{[m_{i+1},h]}\Bigl).\end{multline}
Concerning the first right-hand side term, for every $p_{i+1}>\frac{p_i p}{p-p_i},$ we have
\begin{multline}\label{M1}
\|\nabla v\|_{\dot{\B}_{p,1}^{\frac dp-\frac {d}{p_{i+1}^*}}}^{[\ell,m_{i+1}]}\lesssim
\| v\|_{\dot{\B}_{p,1}^{\frac dp-\frac {d}{p_{i+1}^*}+1}}^{\ell}+\sum^{R}_{\iota=i+1}\| v\|_{\dot{\B}_{p_\iota,1}^{\frac {d}{p_\iota}-\frac {d}{p_{i+1}^*}+1}}^{m_\iota}
\lesssim\varepsilon^{\frac {d}{p_{i+1}^*}-1}\Big(\| v\|_{\dot{\B}_{p,1}^{\frac dp}}^{\ell}+\sum^{R}_{\iota=i+1}\| v\|_{\dot{\B}_{p_\iota,1}^{\frac {d}{p_\iota}}}^{m_\iota}\Big),
\end{multline}
and
\begin{eqnarray}\label{M2}\|\nabla c\|_{\dot{\B}_{p_{i+1},1}^{\frac{d}{p_{i}}-1}}^{m_{i+1}}\lesssim\varepsilon^{2-\frac{d}{p_{i}}+\frac{d}{p_{i+1}}}\| c\|_{\dot{\B}_{p_{i+1},1}^{\frac{d}{p_{i+1}}+2}}^{m_{i+1}},\end{eqnarray}
which implies, noticing that $\frac {d}{p_{i+1}^*}=\frac{d}{p_{i}}-\frac{d}{p_{i+1}}$, 
\begin{multline}\label{mmy1}
\int_0^t \|\nabla v\|_{\dot{\B}_{p,1}^{\frac dp-\frac {d}{p_{i+1}^*}}}^{[\ell,m_{i+1}]}\|\nabla c\|_{\dot{\B}_{p_{i+1},1}^{\frac{d}{p_{i}}-1}}^{m_{i+1}}
\lesssim\varepsilon\big(\| v\|_{L^\infty_T(\dot{\B}_{p,1}^{\frac dp})}^{\ell}+\sum^{R}_{\iota=i+1}\| v\|_{L^\infty_T(\dot{\B}_{p_{\iota},1}^{\frac {d}{p_{\iota}}})}^{m_{\iota}}\big)\| c\|_{L^1_T(\dot{\B}_{p_{i+1},1}^{\frac{d}{p_{i+1}}+2})}^{m_{i+1}}\leq X^2(t).
\end{multline}
Similarly, in light of (\ref{v infty}), we have
\begin{multline}
\int_0^t\|\nabla v\|_{L^\infty}\|\nabla c\|_{\dot{\B}_{p_i,1}^{\frac{d}{p_i}-1}}^{[m_{i},h]}\lesssim\|\nabla v\|_{L^1_T(L^\infty)}\big(\sum^{R}_{\iota=i+1}\| c\|_{L^\infty_T(\dot{\B}_{p_{\iota},1}^{\frac {d}{p_{\iota}}})}^{m_{\iota}}+\varepsilon\|c\|_{L^\infty_T(\dot{\B}_{p,1}^{\frac d2+1})}^{h}\big)\lesssim X^2(t),
\end{multline}
\begin{multline}
\int_0^t \|\nabla c\|_{\dot{\B}_{p,1}^{\frac dp-\frac {d}{p_{i+1}^*}}}^{[\ell,m_{i+1}]}\| v\|_{\dot{\B}_{p_{i+1},1}^{\frac{d}{p_{i}}}}^{m_{i+1}}\lesssim\| c\|_{L^\infty_T(\dot{\B}_{p,1}^{\frac dp-\frac {d}{p_1^*}+1})}^{[\ell,m_{i+1}]}
\| v\|_{L^1_T(\dot{\B}_{p_{i+1},1}^{\frac{d}{p_{i+1}}})}^{m_{i+1}}\\ \lesssim\big(\|c\|_{L^\infty_T(\dot{\B}_{p,1}^{\frac dp})}^{\ell}+\sum^{R}_{\iota=i+1}\| c\|_{L^\infty_T(\dot{\B}_{p_{\iota},1}^{\frac {d}{p_{\iota}}})}^{m_{\iota}}\big)
\| v\|_{L^1_T(\dot{\B}_{p_{i+1},1}^{\frac{d}{p_{i+1}}+1})}^{m_{i+1}}\lesssim X^2(t),
\end{multline}
\begin{eqnarray}
\int_0^t \|\nabla c\|_{L^\infty}\|v\|_{\dot{\B}_{p_i,1}^{\frac{d}{p_i}}}^{[m_{i},h]}\lesssim\|\nabla c\|_{L^2_T(L^\infty)}\big(\sum^{R}_{\iota=i+1}\| v\|_{L^2_T(\dot{\B}_{p_{\iota},1}^{\frac {d}{p_{\iota}}})}^{m_{\iota}}+\varepsilon\|v\|_{L^2_T(\dot{\B}_{2,1}^{\frac d2+1})}^{\ell}\big)\lesssim X^2(t),
\end{eqnarray}
and
\begin{multline}\label{mmy2}
\int_0^t \|\nabla c\|_{\dot{\B}_{p_{i+1},1}^{\frac{d}{p_{i+1}}}}^{[m_{i+1},h]}\|v\|_{\dot{\B}_{p_{i+1},1}^{\frac{d}{p_{i+1}}}}^{[m_{i+1},h]}
\lesssim\big(\| c\|_{L^2_T(\dot{\B}_{2,1}^{\frac{d}{2}+1})}^{h}+\sum^{i+1}_{\iota=1}
\| c\|_{L^2_T(\dot{\B}_{p_{i+1},1}^{\frac{d}{p_{i+1}}+1})}^{m_{\iota}}\big)
\\ \cdot\big(\| v\|_{L^2_T(\dot{\B}_{2,1}^{\frac{d}{2}})}^{h}+\sum^{i+1}_{\iota=1}\| v\|_{L^2_T(\dot{\B}_{p_{i+1},1}^{\frac{d}{p_{i+1}}})}^{m_{\iota}}\big)\lesssim X^2(t).
\end{multline}
Hence we reach
$$\int_0^t\sum_{j\in J^\varepsilon_i}2^{\frac{d}{{p_i}}j}\|\varepsilon\nabla R_j^1\|_{L^{p_i}}\lesssim X^2(t).$$
We can handle the commutator $R^2$ similarly, we have
$$\displaylines{
\sum_{j\in J_i}2^{\frac{d}{p_i}j}\|\varepsilon \nabla R_j^2\|_{L^{p_i}}\leq C\Bigl(\|\nabla c\|_{\dot{\B}_{p,1}^{\frac dp-\frac {d}{p_{i+1}^*}}}^{[\ell,m_{i+1}]}\|\nabla v\|_{\dot{\B}_{p_{i+1},1}^{\frac{d}{p_{i}}-1}}^{m_{i+1}}
+\|\nabla c\|_{L^\infty}\|\nabla v\|_{\dot{\B}_{p_i,1}^{\frac{d}{p_i}-1}}^{[m_{i},h]}
+
\|\nabla v\|_{\dot{\B}_{p,1}^{\frac dp-\frac {d}{p_{i+1}^*}}}^{[\ell,m_{i+1}]}\|c\|_{\dot{\B}_{p_{i+1},1}^{\frac{d}{p_i}}}^{m_{i+1}}\hfill\cr\hfill
+\|\nabla v\|_{L^\infty}\|c\|_{\dot{\B}_{p_i,1}^{\frac{d}{p_i}}}^{[m_{i},h]}+ \|\nabla v\|_{\dot{\B}_{p_{i+1},1}^{\frac{d}{p_{i+1}}}}^{[m_{i+1},h]}\|c\|_{\dot{\B}_{p_{i+1},1}^{\frac{d}{p_{i+1}}}}^{[m_{i+1},h]}\Bigl).}$$
This leads to
$$\int_0^t\sum_{j\in J^\varepsilon_i}2^{\frac{d}{{p_i}}j}\|\varepsilon\nabla R_j^2\|_{L^{p_i}}\lesssim X^2(t).$$
For $R^3_{j}$, Lemma \ref{commutator} yields
$${\sum_{j\in J_i}2^{\frac{d}{p_i}j}\|R_j^3\|_{L^{p_i}}\leq C\Bigl(\|\nabla v\|_{\dot{\B}_{p,1}^{\frac dp-\frac {d}{p_{i+1}^*}}}^{[\ell,m_{i+1}]}\|\nabla v\|_{\dot{\B}_{p_{i+1},1}^{\frac{d}{p_{i}}-1}}^{m_{i+1}}+\|\nabla v\|_{L^\infty}\|\nabla v\|_{\dot{\B}_{p_i,1}^{\frac{d}{p_i}-1}}^{[m_{i},h]}
+ \|\nabla v\|_{\dot{\B}_{p_{i+1},1}^{\frac{d}{p_{i+1}}}}^{[m_{i+1},h]}\|v\|_{\dot{\B}_{p_{i+1},1}^{\frac{d}{p_{i+1}}}}^{[m_{i+1},h]}\Bigl).}$$
It holds
\begin{multline*}
\varepsilon\int_0^t \|\nabla v\|_{\dot{\B}_{p,1}^{\frac dp-\frac {d}{p_{i+1}^*}}}^{[\ell,m_{i+1}]}\|\nabla v\|_{\dot{\B}_{p_{i+1},1}^{\frac{d}{p_{i}}-1}}^{m_{i+1}}
\lesssim\varepsilon\big(\| v\|_{L^\infty_T(\dot{\B}_{p,1}^{\frac dp})}^{\ell}+\sum^{R}_{\iota=i+1}\| v\|_{L^\infty_T(\dot{\B}_{p_{\iota},1}^{\frac {d}{p_{\iota}}})}^{m_{\iota}}\big)\| v\|_{L^1_T(\dot{\B}_{p_{i+1},1}^{\frac{d}{p_{i+1}}+1})}^{m_{i+1}}\lesssim X^2(t),
\end{multline*}
\begin{eqnarray*}
\varepsilon\int_0^t \|\nabla v\|_{L^\infty}\|v\|_{\dot{\B}_{p_i,1}^{\frac{d}{p_i}}}^{[m_{i},h]}\lesssim\|\nabla v\|_{L^1_T(L^\infty)}\big(\sum^{R}_{\iota=i+1}\varepsilon\| v\|_{L^\infty_T(\dot{\B}_{p_{\iota},1}^{\frac {d}{p_{\iota}}})}^{m_{\iota}}+\varepsilon^{2}\|v\|_{L^\infty_T(\dot{\B}_{2,1}^{\frac d2+1})}^{h}\big)\lesssim X^2(t),
\end{eqnarray*}
\begin{multline*}
\varepsilon\int_0^t \|\nabla v\|_{\dot{\B}_{p_{i+1},1}^{\frac{d}{p_{i+1}}}}^{[m_{i+1},h]}\|v\|_{\dot{\B}_{p_{i+1},1}^{\frac{d}{p_{i+1}}}}^{[m_{i+1},h]}
\lesssim\big(\varepsilon\|v\|_{L^2_T(\dot{\B}_{2,1}^{\frac{d}{2}+1})}^{h}+\sum^{i+1}_{\iota=1}\varepsilon
\|v\|_{L^2_T(\dot{\B}_{p_{i+1},1}^{\frac{d}{p_{i+1}}+1})}^{m_{\iota}}\big)
\\ \cdot\big(\| v\|_{L^2_T(\dot{\B}_{2,1}^{\frac{d}{2}})}^{h}+\sum^{i+1}_{\iota=1}\| v\|_{L^2_T(\dot{\B}_{p_{i+1},1}^{\frac{d}{p_{i+1}}})}^{m_{\iota}}\big)\lesssim X^2(t).
\end{multline*}
Gathering the estimates, we arrive at
$$\varepsilon\int_0^t\sum_{j\in J^\varepsilon_i}2^{\frac{d}{{p_i}}j}\|R_j^3\|_{L^{p_i}}\lesssim X^2(t).$$
As for the term $\partial_t(c\nabla c)$ in $Q^3_j$, we find
\begin{eqnarray}\label{partt}\partial_t(c\nabla c)=-\nabla\left(\check{\gamma}\bar cc\div v+cv\cdot\nabla c+\check{\gamma}c^2\div v\right).\end{eqnarray}
The quadratic term can be handled exactly as before. For the cubic term $cv\cdot c$, we have
$$\ddj(cv\cdot\nabla c)=\dot S_{j-1}(cv)\cdot\nabla\ddj c+\tilde R_j^1,\quad
\tilde R_j^1:=\dot S_{j-1}(cv)\cdot\nabla\ddj c-\ddj(cv\cdot\nabla c).$$
Consequently computations analoguous to (\ref{SSDDFF}) and (\ref{SSDDFFGG}) show that
\begin{eqnarray*}
\int_0^t\sum_{j\in J^\varepsilon_i}2^{\frac{d}{{p_i}}j}\|\varepsilon \nabla(\dot{S}_{j-1}(cv)\cdot\nabla \ddj c)\|_{L^{p_i}}&\lesssim&\int_0^t\sup_{j\in J^\varepsilon_i}\|\dot{S}_{j-1}(cv)\|_{L^\infty}\|c\|_{\dot{\B}_{p_i,1}^{\frac{d}{p_i}+1}}^{m_i}\\
&\leq&\|c\|_{L^\infty(L^\infty)}\|v\|_{L^2(L^\infty)}\|c\|_{L^2_T(\dot{\B}_{p_i,1}^{\frac{d}{p_i}+1})}^{m_i}\lesssim X^3(t),
\end{eqnarray*}
while for the commutator term we have
$$\displaylines{
\sum_{j\in J_i}2^{\frac{d}{p_i}j}\|\varepsilon \nabla \tilde R_j^1\|_{L^{p_i}}\leq C\Bigl(\|\nabla (cv)\|_{\dot{\B}_{p,1}^{\frac dp-\frac {d}{p_{i+1}^*}}}^{[\ell,m_{i+1}]}\|\nabla c\|_{\dot{\B}_{p_{i+1},1}^{\frac{d}{p_{i}}-1}}^{m_{i+1}}
+\|\nabla (cv)\|_{L^\infty}\|\nabla c\|_{\dot{\B}_{p_i,1}^{\frac{d}{p_i}-1}}^{[m_{i},h]}
\hfill\cr\hfill+
\|\nabla c\|_{\dot{\B}_{p,1}^{\frac dp-\frac {d}{p_{i+1}^*}}}^{[\ell,m_{i+1}]}\|(cv)\|_{\dot{\B}_{p_{i+1},1}^{\frac{d}{p_i}}}^{m_{i+1}}
+\|\nabla c\|_{L^\infty}\|(cv)\|_{\dot{\B}_{p_i,1}^{\frac{d}{p_i}}}^{[m_{i},h]}+ \|\nabla c\|_{\dot{\B}_{p_{i+1},1}^{\frac{d}{p_{i+1}}}}^{[m_{i+1},h]}\|(cv)\|_{\dot{\B}_{p_{i+1},1}^{\frac{d}{p_{i+1}}}}^{[m_{i+1},h]}\Bigl).}$$
We only treat the first two right-hand side terms, the other being easier to handle.  Concerning the first term, for every $p_{i+1}>\frac{p_i p}{p-p_i},$ using (\ref{M1}) yields
\begin{multline*}
\|\nabla (cv)\|_{\dot{\B}_{p,1}^{\frac dp-\frac {d}{p_{i+1}^*}}}^{[\ell,m_{i+1}]}\lesssim
\sum^{R}_{\iota=i+1}\| cv\|_{\dot{\B}_{p_\iota,1}^{\frac {d}{p_\iota}-\frac {d}{p_{i+1}^*}+1}}^{\ell}+\| cv\|_{\dot{\B}_{p,1}^{\frac dp-\frac {d}{p_{i+1}^*}+1}}^{m_\iota}\\
\lesssim\varepsilon^{\frac {d}{p_{i+1}^*}-1}\| cv\|_{\dot{\B}_{p,1}^{\frac dp}}^{\ell}+\varepsilon^{\frac {d}{p_{i+1}^*}}\sum^{R}_{\iota=i+1}\| cv\|_{\dot{\B}_{p_\iota,1}^{\frac {d}{p_\iota}+1}}^{m_\iota}.
\end{multline*}
We have
$$\| cv\|_{\dot{\B}_{p,1}^{\frac dp}}^{\ell}\lesssim\|c\|_{\dot{\B}_{p,1}^{\frac dp}}\|v\|_{\dot{\B}_{p,1}^{\frac dp}}$$
and, recalling (\ref{M2}), we obtain 
\begin{multline}
\int_0^t \varepsilon^{\frac {d}{p_{i+1}^*}-1}\| cv\|_{\dot{\B}_{p,1}^{\frac dp}}^{\ell}\|\nabla c\|_{\dot{\B}_{p_{i+1},1}^{\frac{d}{p_{i}}-1}}^{m_{i+1}}
\lesssim\| c\|_{L^\infty_T(\dot{\B}_{p,1}^{\frac dp})}\| v\|_{L^2_T(\dot{\B}_{p,1}^{\frac dp})}\| c\|_{L^2_T(\dot{\B}_{p_{i+1},1}^{\frac{d}{p_{i+1}}+1})}^{m_{i+1}}\leq X^3(t).
\end{multline}
On the other hand, noticing that
$$\sum^{R}_{\iota=i+1}\| cv\|_{\dot{\B}_{p_\iota,1}^{\frac {d}{p_\iota}+1}}^{m_\iota}\lesssim\sum^{R}_{\iota=i+1}\| c\div v\|_{\dot{\B}_{p_\iota,1}^{\frac {d}{p_\iota}}}^{m_\iota}+\sum^{R}_{\iota=i+1}\| v\cdot\nabla c\|_{\dot{\B}_{p_\iota,1}^{\frac {d}{p_\iota}}}^{m_\iota},$$
by repeating the computations carried out in \eqref{SSDDFF} and \eqref{SSDDFFGG}–\eqref{mmy2}, we obtain
\begin{multline}
\int_0^t \varepsilon^{\frac {d}{p_{i+1}^*}}\sum^{R}_{\iota=i+1}\| cv\|_{\dot{\B}_{p_\iota,1}^{\frac {d}{p_\iota}+1}}^{m_\iota}\|\nabla c\|_{\dot{\B}_{p_{i+1},1}^{\frac{d}{p_{i}}-1}}^{m_{i+1}}
\lesssim\sum^{R}_{\iota=i+1}\| cv\|_{L^1_T(\dot{\B}_{p_\iota,1}^{\frac {d}{p_\iota}+1})}^{m_\iota}\|c\|_{L^\infty_T(\dot{\B}_{p_{i+1},1}^{\frac{d}{p_{i+1}}})}^{m_{i+1}}\leq X^3(t).
\end{multline}
Combining the above inequalities, we obtain
$$\displaylines{
\sum_{j\in J_i}2^{\frac{d}{p_i}j}\|\varepsilon \nabla \tilde R_j^1\|_{L^{p_i}}\leq C X^3(t)\quad\Longrightarrow \int_0^t\sum_{j\in J^\varepsilon_i}2^{\frac{d}{{p_i}}j}\|\varepsilon \ddj\nabla(cv\cdot\nabla  c)\|_{L^{p_i}}\leq C X^3(t).}$$
The other cubic terms in (\ref{partt}) can be bounded in a similar fashion. This completes estimates for $Q^3_j$.

{\bf Estimate for $Q^4_{j}$}. It remains to estimate the term $\div(c\cdot\nabla c)$. We have
% Recall that
% \begin{eqnarray}
% c\div c=\dot S_{j-1}c\cdot\nabla\ddj c+R_j^4;
% \end{eqnarray}
\begin{eqnarray}
\ddj(\div(c\cdot\nabla c))=\ddj(\nabla c\cdot \nabla c)+\dot S_{j-1}c\cdot\Delta\ddj c+\tilde R_j^4,\:\:\tilde R_j^4
:=\dot S_{j-1}c\Delta\ddj c-\ddj(c\Delta c).\end{eqnarray}
Then, the product estimates in Lemma \ref{product} allow us to get
\begin{multline*}\sum_{j\in J_i}2^{\frac{d}{p_i}j}\|\div (\nabla c\cdot \nabla c)\|_{L^{p_i}}\leq C\Bigl(\|\nabla c\|_{\dot{\B}_{p,1}^{\frac dp-\frac {d}{p_{i+1}^*}}}^{[\ell,m_{i+1}]}\|\nabla c\|_{\dot{\B}_{p_{i+1},1}^{\frac{d}{p_{i}}}}^{m_{i+1}}+\|\nabla c\|_{L^\infty}\|\nabla c\|_{\dot{\B}_{p_i,1}^{\frac{d}{p_i}}}^{[m_{i},h]}\\
+ \|\nabla c\|_{\dot{\B}_{p_{i+1},1}^{\frac{d}{p_{i+1}}}}^{[m_{i+1},h]}\|c\|_{\dot{\B}_{p_{i+1},1}^{\frac{d}{p_{i+1}}+1}}^{[m_{i+1},h]}\Bigl),\end{multline*}
and
\begin{multline*}
\int_0^t\sum_{j\in J_i}2^{\frac{d}{{p_i}}j}\|\div(\dot{S}_{j-1}c\nabla \ddj c)\|_{L^{p_i}}\lesssim\int_0^t\sup_{j\in J^\varepsilon_i}\|\dot{S}_{j-1}c\|_{L^\infty}\|c\|_{\dot{\B}_{p_i,1}^{\frac{d}{p_i}+2}}^{m_i}\\
\leq\|c\|_{L^\infty_T( L^\infty)}\|c\|_{L^1_T(\dot{\B}_{p_i,1}^{\frac{d}{p_i}+2})}^{m_i}\lesssim X^2(t).
\end{multline*}
For $R^4_{j}$, there holds
$$\displaylines{\sum_{j\in J_i}2^{\frac{d}{p_i}j}\|\div R_j^4\|_{L^{p_i}}\leq C\Bigl(\|\nabla c\|_{\dot{\B}_{p,1}^{\frac dp-\frac {d}{p_{i+1}^*}}}^{[\ell,m_{i+1}]}\|\nabla c\|_{\dot{\B}_{p_{i+1},1}^{\frac{d}{p_{i}}}}^{m_{i+1}}+\|\nabla c\|_{L^\infty}\|\nabla c\|_{\dot{\B}_{p_i,1}^{\frac{d}{p_i}}}^{[m_{i},h]}
+ \|\nabla c\|_{\dot{\B}_{p_{i+1},1}^{\frac{d}{p_{i+1}}}}^{[m_{i+1},h]}\|c\|_{\dot{\B}_{p_{i+1},1}^{\frac{d}{p_{i+1}}+1}}^{[m_{i+1},h]}\Bigl).}$$
Then it holds that
\begin{multline*}
\int_0^t \|\nabla c\|_{\dot{\B}_{p,1}^{\frac dp-\frac {d}{p_{i+1}^*}}}^{[\ell,m_{i+1}]}\|\nabla c\|_{\dot{\B}_{p_{i+1},1}^{\frac{d}{p_{i}}}}^{m_{i+1}}
\lesssim\big(\|c\|_{L^\infty_T(\dot{\B}_{p,1}^{\frac dp})}^{\ell}+\sum^{R}_{\iota=i+1}\| c\|_{L^\infty_T(\dot{\B}_{p_{\iota},1}^{\frac {d}{p_{\iota}}})}^{m_{\iota}}\big)\frac{1}{\varepsilon}\| c\|_{L^1_T(\dot{\B}_{p_{i+1},1}^{\frac{d}{p_{i+1}}+2})}^{m_{i+1}}\lesssim X^2(t),
\end{multline*}
\begin{eqnarray*}
\int_0^t \|\nabla c\|_{L^\infty}\|c\|_{\dot{\B}_{p_i,1}^{\frac{d}{p_i}+1}}^{[m_{i},h]}\lesssim\|\nabla c\|_{L^2 L^\infty}\big(\sum^{R}_{\iota=i+1}\|c\|_{L^2_T(\dot{\B}_{p_{\iota},1}^{\frac {d}{p_{\iota}+2}})}^{m_{\iota}}+\|c\|_{L^2_T(\dot{\B}_{2,1}^{\frac d2+1})}^{\ell}\big)\lesssim X^2(t),
\end{eqnarray*}
\begin{eqnarray*}
\int_0^t \|\nabla c\|_{\dot{\B}_{p_{i+1},1}^{\frac{d}{p_{i+1}}}}^{[m_{i+1},h]}\|c\|_{\dot{\B}_{p_{i+1},1}^{\frac{d}{p_{i+1}}+1}}^{[m_{i+1},h]}
\lesssim\big(\|c\|_{L^2_T(\dot{\B}_{2,1}^{\frac{d}{2}+1})}^{h}+\sum^{i+1}_{\iota=1}
\|c\|_{L^2_T(\dot{\B}_{p_{i+1},1}^{\frac{d}{p_{i+1}}+1})}^{m_{\iota}}\big)^2
\lesssim X^2(t),
\end{eqnarray*}
which gives rise to
$$\int_0^t\sum_{j\in J^\varepsilon_i}2^{\frac{d}{{p_i}}j}\|\div R_j^4\|_{L^{p_i}}\lesssim X^2(t).$$
The above computations yield
\begin{eqnarray}\label{priori2}
X^{m_i}(t) \lesssim X_0+X^2(t)+X^3(t).
\end{eqnarray}

\subsubsection{Low frequencies estimates in \texorpdfstring{$L^p$}{Lp}}
Estimating the linear part is standard, see \cite{CBD3}. Using the damping mode $w$, we obtain 
\begin{eqnarray}
X^\ell(t)\leq
X_{0}+\int_0^t\sum_{j\leq \J-N_0 R}2^{\frac{d}{{p}}j}\| (Q^3_{j},Q^4_{j})\|_{L^{p}}.
\end{eqnarray}
The nonlinear analysis in the low-frequency regime is much simpler as we do not need to perform commutators estimates nor use Hölder inequality.
Indeed, one just rely on product law and Sobolev embeddings.  For instance, one has
\begin{eqnarray*}
\int_0^t\sum_{j< \J-N_0 R}2^{\frac{d}{p}j}\|\ddj(v\cdot\nabla c)\|_{L^p}
\lesssim\|v\|_{L^2_T(\dot{\B}^{\frac{d}{p}}_{p,1})}\| c\|_{L^2_T(\dot{\B}^{\frac{d}{p}+1}_{p,1})}\lesssim X^2(t)
\end{eqnarray*}
and
\begin{eqnarray*}
\int_0^t\sum_{j< \J-N_0 R}2^{\frac{d}{p}j}\|\ddj\div(c\cdot\nabla c)\|_{L^p}
\lesssim\|c\|_{L^\infty_T(\dot{\B}^{\frac{d}{p}}_{p,1})}\| c\|_{L^1_T(\dot{\B}^{\frac{d}{p}+2}_{p,1})}\lesssim X^2(t)
\end{eqnarray*}
where we used the embedding $\dot{\B}^{s}_{2,1}\hookrightarrow \dot{\B}^{s-\frac d2+\frac dp}_{p,1}$. The other terms can be controlled in a similar fashion. We obtain
\begin{eqnarray}\label{priori3}
X^\ell(t)\lesssim X_0+X^2(t)+X^3(t).
\end{eqnarray}
Combining (\ref{priori1}), (\ref{priori2}) and (\ref{priori3}), we finish the proof of Lemma \ref{priori}.

\subsection{Global well-posedness}
In this subsection, we outline the main arguments underlying the global-in-time existence and uniqueness of solutions to system \eqref{CED4}.

From the a priori estimate in Lemma \ref{priori}, the  existence of global-in-time solutions can be established by using a standard bootstrap argument relying on the smallness of the initial data. We refer to \cite[Section 3.2]{CBD3} for further details.

Concerning the uniqueness of the solution, it is sufficient to consider the case $\varepsilon=1$. We define $(\delta c,\delta v)=(c_1-c_2, v_1-v_2)$ where $(c_1,v_2)$ and $(c_2,v_2)$ are global-in-time solutions, associated to the initial data $(c_{1,0},v_{1,0})$ and $(c_{2,0},v_{2,0})$, respectively, verifying the regularity properties from Theorem \ref{thm1}. We have
\begin{equation}\label{delta z}
\left\{
\begin{array}{l}\partial_t{\delta c}+\check{\gamma}\bar c{\div}\delta v=\delta Q^1\triangleq Q^1(c_1,v_1)-Q^1(c_2,v_2), \\ [1mm]
\partial_t\delta v+\check{\gamma}\bar c\nabla{\delta c}+\delta v=\delta Q^2\triangleq Q^2(c_1,v_1)-Q^2(c_2,v_2).\\[1mm]
 \end{array} \right.
\end{equation}
To establish uniqueness, we aim to show that $\delta c = \delta v = 0$. 
This is achieved by estimating \eqref{delta z}, using the properties of 
the solutions and applying Gronwall’s lemma.
We denote
$$\delta X(t)\triangleq\|(\delta c,\delta  v)\|_{\dot{\B}_{2,1}^{\frac{d}{2}}}^h+\sum^{R}_{i=1}\|(\delta c,\delta  v)\|_{\dot{\B}_{{p_i},1}^{\frac{d}{{p_i}}}}^{m_{i} }+\|(\delta c,\delta  v)\|_{\dot{\B}_{p,1}^{\frac{d}{p}}}^{\ell}$$
and
$$\delta X_{0}\triangleq\|(\delta c_{0},\delta  v_{0})\|_{\dot{\B}_{2,1}^{\frac{d}{2}}}^h+\sum^{R}_{i=1}\|(\delta c_{0},\delta  v_{0})\|_{\dot{\B}_{{p_i},1}^{\frac{d}{{p_i}}}}^{m_{i} }+\|(\delta c_{0},\delta  v_{0})\|_{\dot{\B}_{p,1}^{\frac{d}{p}}}^{\ell}.$$
We have the following result.
\begin{lem}\label{Error}
Let $d\geq1$, $p\in(2,\infty)$ and $(p_i)\in\mathcal{S}_{p_i,R}$. Let $(c_1,v_1)$ and $(c_2,v_2)$ be two solutions of (\ref{CED4}) associated to the initial data $(c_{1,0},v_{1,0})$ and $(c_{2,0},v_{2,0})$, respectively. The following estimate holds true
\begin{eqnarray}\label{error}
\delta X(t)\lesssim\delta X_{0}+\int_0^t \tilde{X}(s)\delta X(s)
\end{eqnarray}
where
\begin{eqnarray*}\tilde{X}(s)\triangleq\|(c_{i}, v_{i})\|_{\dot{\B}_{2,1}^{\frac{d}{2}+1}}^h+\sum^{R}_{i=1}\|( c_{i},  v_{i})\|_{\dot{\B}_{{p_i},1}^{\frac{d}{{p_i}}}}+\|( c_{i}, v_{i})\|_{\dot{\B}_{p,1}^{\frac{d}{p}}}^{\ell}.\end{eqnarray*}
\end{lem}
Employing Lemma \ref{Error}, we are able to conclude the uniqueness using Gronwall's inequality and the fact that
$$\sup_{t\in(0,T]}\int^t_0\tilde X(s)\leq C$$
for $T\leq c$ with $c$ sufficient small. See \cite[Section 3.3]{CBD3} for additional details.

In the next sections, we provide the proof of Lemma \ref{Error}.
\subsubsection{Error estimates in high frequencies}

The error system reads
\begin{equation*}
\left\{
\begin{array}{l}\partial_t{\delta c}_{j}+\dot S_{j-1}v_{1}\cdot\nabla{\delta c}_{j}+{\div}\delta v_{j}=\delta R_j^1+\delta R_j^2-\dot S_{j-1}c_{1}\div{\delta v}_{j}+\delta \tilde Q^1_{j}, \\ [1mm]
\partial_t\delta v_{j}+\dot S_{j-1}v_{1}\cdot\nabla \delta v_{j}+\nabla{\delta c}_{j}+\delta v_{j}=\delta R_j^3+\delta R_j^4-\dot S_{j-1}c_{1}\nabla \delta c_{j}+\delta \tilde Q^2_{j},\\[1mm]
 \end{array} \right.
\end{equation*}
where
\begin{eqnarray*}\delta R_j^1=\dot S_{j-1}v_{1}\cdot\nabla\ddj\delta  c-\ddj(v_{1}\cdot\nabla \delta c), \delta R_j^2=\dot S_{j-1}c_{1}\div\ddj \delta v-\ddj(c_{1}\div \delta v),\end{eqnarray*}
\begin{eqnarray*}\delta R_j^3=\varepsilon^2\big(\dot S_{j-1}v_{1}\cdot\nabla\ddj\delta  v-\ddj(v_{1}\cdot\nabla\delta  v)\big), \delta R_j^4=\dot S_{j-1}c_{1}\nabla\ddj\delta  c-\ddj(c_{1}\nabla \delta c)\end{eqnarray*}
and
\begin{eqnarray*}\delta \tilde Q^1_j=\ddj\big(\delta v\cdot\nabla c_2-\delta c\cdot\nabla v_2\big),\: \delta \tilde Q^2_j=\ddj\big(\delta c\cdot\nabla c_2-\delta v\cdot\nabla v_2\big).\end{eqnarray*}
One has
\begin{multline}\label{high f1}
\|(\delta c,\delta  v)(t)\|_{\dot{\B}_{2,1}^{\frac{d}{2}}}^h\lesssim
\delta X_{0}+\int_0^t\|\nabla v_{1}\|_{L^\infty}\|(\delta c, \delta v)\|_{\dot{\B}_{2,1}^{\frac{d}{2}}}^h\\
+\int_0^t\sum_{j\geq \J}2^{\frac{d}{2}j}\left\|(\delta R_j^1,\delta R_j^2,\delta \tilde Q^1_{j},\delta R_j^3,\delta R_j^4,\delta \tilde Q^2_{j})\right\|_{L^2}\\+\int_0^t\sum_{j\geq \J}2^{\frac{d}{2}j}\left\|( \delta Q^1_{j},\delta Q^2_{j})\right\|_{L^2}.
\end{multline}
We have
\begin{eqnarray}\label{luq}
\int_0^t\|\nabla v_{1}\|_{L^\infty}\|(\delta c, \delta v)\|_{\dot{\B}_{2,1}^{\frac{d}{2}}}^h
\lesssim\int_0^t\tilde{X}(s)\delta X(s)
\end{eqnarray}
where we used that
$$\|\nabla v_{1}\|_{L^\infty}\lesssim\|v_{i}\|_{\dot{\B}_{2,1}^{\frac{d}{2}+1}}^h+\sum^{R}_{i=1}\| v_{i}\|_{\dot{\B}_{{p_i},1}^{\frac{d}{{p_i}}+1}}^{m_{i}}+\|v_{i}\|_{\dot{\B}_{p,1}^{\frac{d}{p}+1}}^{\ell}\lesssim\tilde{X}(s).$$
Then it is sufficient to pay attention to $\delta R_j^1$ and $\delta \tilde Q^1_{j}$. Lemma \ref{commutator} gives
\begin{multline}\label{uniqee}
\sum_{j\geq \J}2^{\frac{d}{2}j}\|\delta R_j^1\|_{L^2}\leq C\Bigl(\|\nabla v_{1}\|_{\dot{\B}_{p,1}^{\frac dp-\frac {d}{p_1^*}}}^{[\ell,m_1]}\|\nabla\delta c\|_{\dot{\B}_{p_1,1}^{\frac{d}{2}-1}}^{m_1}
+\|\nabla v_{1}\|_{L^\infty}\|\nabla\delta  c\|_{\dot{\B}_{2,1}^{\frac{d}{2}-1}}^{h}\hfill\cr\hfill+
\|\nabla\delta  c\|_{\dot{\B}_{p,1}^{\frac dp-\frac {d}{p_1^*}}}^{[\ell,m_1]}\|v_{1}\|_{\dot{\B}_{p_1,1}^{\frac{d}{2}}}^{m_1}
+\|\nabla\delta  c\|_{\dot{\B}_{\infty,1}^{-1}}\|v_{1}\|_{\dot{\B}_{2,1}^{\frac{d}{2}}}^{h}+ \|\nabla\delta  c\|_{\dot{\B}_{p_1,1}^{\frac{d}{p_1}-1}}^{[m_1,h]}\|v_{1}\|_{\dot{\B}_{p_1,1}^{\frac{d}{p_1}+1}}^{[m_1,h]}\Bigl).
\end{multline}
We just deal with the first right-hand side term, the other terms can be handled as in the previous sections. Recalling that $r=1-\frac{d}{2}+\frac{d}{p_1}$, it holds
$$\|\nabla v_{1}\|_{\dot{\B}_{p,1}^{\frac dp-\frac {d}{p_1^*}}}^{[\ell,m_1]}\lesssim \| v_{1}\|_{\dot{\B}_{p,1}^{\frac dp}}^{\ell}+\sum^{R}_{i=1}\| v_{1}\|_{\dot{\B}_{p_i,1}^{\frac {d}{p_i}}}^{m_i}\lesssim\tilde X(s),$$
$$\|\nabla\delta c\|_{\dot{\B}_{p_1,1}^{\frac{d}{2}-1}}^{m_1}\lesssim\|\delta c\|_{\dot{\B}_{p_1,1}^{\frac{d}{p_1}}}^{m_1}\lesssim\delta X(s),$$
which gives
\begin{eqnarray}\label{luq2}
\int_0^t\|\nabla v_{1}\|_{\dot{\B}_{p,1}^{\frac dp-\frac {d}{p_1^*}}}^{[\ell,m_1]}\|\nabla\delta c\|_{\dot{\B}_{p_1,1}^{\frac{d}{2}-1}}^{m_1}
\lesssim\int_0^t\tilde{X}(s)\delta X(s).
\end{eqnarray}
By treating the other terms in the same manner as in (\ref{luq}) and (\ref{luq2}), we have
\begin{eqnarray*}
\int_0^t\sum_{j\geq \J}2^{\frac{d}{2}j}\|\delta R_j^1\|_{L^2}\leq \int_0^t \tilde{X}(s)\delta X(s).
\end{eqnarray*}
For $\delta \tilde Q^1_{j}$, we focus on $\delta v\cdot\nabla c_2$. We decompose
\begin{eqnarray}\label{S1}\delta v\cdot\nabla c_2=\delta v^h\cdot\nabla c_2+(\sum_{\iota=1}^R\delta v^{m_\iota}+\delta v^\ell)\cdot\nabla c_2.\end{eqnarray}
For the high-frequency term, we have
\begin{eqnarray*}
\sum_{j\geq \J}2^{\frac{d}{2}j}\|\ddj(\delta v^h\cdot\nabla c_2)\|_{L^2}\leq C\|\delta v^h\|_{\dot{\B}_{2,1}^{\frac d2}}\|\nabla c_{2}\|_{\dot{\B}_{\infty,1}^{0}}\lesssim\tilde{X}(s)\delta X(s).
\end{eqnarray*}
For the other terms, we decompose
\begin{eqnarray}\label{S2}
\ddj\big(\delta v^{m_\iota}\cdot\nabla c_2\big)=\dot S_{j-1}\delta v^{m_\iota}\cdot\nabla\ddj c_{2}+\delta \bar R_j^1(v^{m_\iota},c_2)\,\,
\end{eqnarray}
and
\begin{eqnarray}\label{S3}
\ddj\big(\delta v^\ell\cdot\nabla c_2\big)=\dot S_{j-1}\delta v^\ell\cdot\nabla\ddj c_{2}+\delta \bar R_j^1(v^\ell,c_2),
\end{eqnarray}
where
$$\delta \bar R_j^1(f,c_2)=\dot S_{j-1}f\cdot\nabla\ddj c_2-\ddj(f\cdot\nabla c_2).$$
We have
\begin{eqnarray*}
\sum_{j\geq \J}2^{\frac{d}{2}j}\|\dot S_{j-1}\delta v^{m_\iota}\cdot\nabla\ddj c_{2}\|_{L^2}\leq C\|\dot S_{j-1}\delta v^{m_\iota}\|_{L^\infty}\|c_2\|_{\dot{\B}_{2,1}^{\frac{d}{2}+1}}^{h}
\lesssim\tilde{X}(s)\delta X(s)
\end{eqnarray*}
and for the remainder, Lemma \ref{commutator} yields
\begin{multline*}
\sum_{j\geq \J}2^{\frac{d}{2}j}\|\delta\bar R_j^1(v^{m_\iota},c_2)\|_{L^2}\leq C\Bigl(\|\nabla\delta v^{m_\iota}\|_{\dot{\B}_{p,1}^{\frac dp-\frac {d}{p_1^*}}}^{[\ell,m_1]}\|\nabla c_{2}\|_{\dot{\B}_{p_1,1}^{\frac{d}{2}-1}}^{m_1}
+\|\nabla\delta v^{m_\iota}\|_{L^\infty}\|\nabla c_{2}\|_{\dot{\B}_{2,1}^{\frac{d}{2}-1}}^{h}\hfill\cr\hfill+
\|\nabla  c_{2}\|_{\dot{\B}_{p,1}^{\frac dp-\frac {d}{p_1^*}}}^{[\ell,m_1]}\|\delta v^{m_\iota}\|_{\dot{\B}_{p_1,1}^{\frac{d}{2}}}^{m_1}
+\|\nabla  c_{2}\|_{\dot{\B}_{\infty,1}^{-1}}\|\delta v^{m_\iota}\|_{\dot{\B}_{2,1}^{\frac{d}{2}}}^{h}+ \|\nabla\delta  v^{m_\iota}\|_{\dot{\B}_{p_1,1}^{\frac{d}{p_1}-1}}^{[m_1,h]}\|c_{2}\|_{\dot{\B}_{p_1,1}^{\frac{d}{p_1}+1}}^{[m_1,h]}\Bigl).
\end{multline*} 
Compared to \eqref{uniqee}, we have to handle the second right-hand side term differently. We have
$$\|\nabla\delta v^{m_\iota}\|_{L^\infty}\lesssim\|\nabla\delta v\|^{m_\iota}_{\dot{\B}_{p_1,1}^{\frac{d}{p_1}}}\lesssim\|\delta v\|^{m_\iota}_{\dot{\B}_{p_1,1}^{\frac{d}{p_1}}};\quad\|\nabla c_{2}\|_{\dot{\B}_{2,1}^{\frac{d}{2}-1}}^{h}\lesssim\|\nabla c_{2}\|_{\dot{\B}_{2,1}^{\frac{d}{2}+1}}^{h},$$
which implies that
\begin{eqnarray*}
\|\nabla\delta v^{m_\iota}\|_{L^\infty}\|\nabla c_{2}\|_{\dot{\B}_{2,1}^{\frac{d}{2}-1}}^{h}
\lesssim\tilde{X}(s)\delta X(s)
\end{eqnarray*}
and gives
\begin{eqnarray*}
\int_0^t\sum_{j\geq \J}2^{\frac{d}{2}j}\|\delta\bar R_j^1\|_{L^2}\leq \int_0^t \tilde{X}(s)\delta X(s).
\end{eqnarray*}
The remaining terms can be treated similarly. Thus we reach
\begin{eqnarray}
\|(\delta c,\delta  v)\|_{\dot{\B}_{2,1}^{\frac{d}{2}}}^h\lesssim\delta X_{0}+\int_0^t \tilde{X}(s)\delta X(s).
\end{eqnarray}

\subsubsection{Error estimates in medium frequencies}
We define the error damped mode $\delta w= \delta v+\check{\gamma}\bar c\nabla \delta c$. The system reads
\begin{equation}
\left\{
\begin{array}{l}\partial_{t}\delta w+ \delta w=-\nabla\div \delta v+\nabla \delta Q^1+\delta Q^2, \\ 
\partial_t\delta  c-\check{\gamma}\bar c\Delta\delta  c=-\textrm{div}\delta w+\delta Q^1.\\[1mm]
 \end{array} \right.
\end{equation}
Repeating the computations from Section~4.1.2, we obtain
\begin{eqnarray}
\| (\delta c,\delta v)\|^{m_i}_{L^\infty_T(\dot{\mathbb{B}}_{p_i,1}^{\frac{d}{p_i}})}\leq\delta X_{0}+\int_0^t\sum_{j\in J_i}2^{\frac{d}{{p_i}}j}\| (\delta Q^1_{j},\delta Q^2_{j})\|_{L^{p_i}}.
\end{eqnarray}
We focus on the term $v_1\cdot\nabla\delta c$ in $\delta Q^1_{j}$. We have
\begin{multline}\label{gubb}
\int_0^t\sum_{j\in J_i}2^{\frac{d}{p_i}j}\|\dot S_{j-1}v_1\cdot\ddj\nabla\delta c\|_{L^{p_i}}\leq C\int_0^t\|\dot S_{j-1}v_1\|_{L^\infty}\|\nabla\delta c\|_{\dot{\B}_{p_i,1}^{\frac{d}{p_i}}}^{m_i}
\leq \int^t_0\delta X(s)\tilde X(s).
\end{multline}
For the commutator,  Lemma \ref{commutator} gives
\begin{multline}
\sum_{j\in J_i}2^{\frac{d}{p_i}j}\|\delta R_j^1\|_{L^{p_i}}\leq C\Bigl(\|\nabla v_1\|_{\dot{\B}_{p,1}^{\frac dp-\frac {d}{p_{i+1}^*}}}^{[\ell,m_{i+1}]}\|\nabla \delta c\|_{\dot{\B}_{p_{i+1},1}^{\frac{d}{p_{i}}-1}}^{m_{i+1}}
+\|\nabla v_1\|_{L^\infty}\|\nabla \delta c\|_{\dot{\B}_{p_i,1}^{\frac{d}{p_i}-1}}^{[m_{i},h]}\hfill\cr\hfill
+
\|\nabla\delta  c\|_{\dot{\B}_{p,1}^{\frac dp-\frac {d}{p_{i+1}^*}}}^{[\ell,m_{i+1}]}\|v_1\|_{\dot{\B}_{p_{i+1},1}^{\frac{d}{p_i}}}^{m_{i+1}}
+\|\nabla \delta c\|_{\dot{\B}_{\infty,1}^{-1}}\|v_1\|_{\dot{\B}_{p_i,1}^{\frac{d}{p_i}}}^{[m_{i},h]}+ \|\nabla\delta  c\|_{\dot{\B}_{p_{i+1},1}^{\frac{d}{p_{i+1}}-1}}^{[m_{i+1},h]}\|v_1\|_{\dot{\B}_{p_{i+1},1}^{\frac{d}{p_{i+1}}+1}}^{[m_{i+1},h]}\Bigl).\end{multline}
For the first term, using (\ref{M1}), we get
\begin{eqnarray}
\|\nabla v_1\|_{\dot{\B}_{p,1}^{\frac dp-\frac {d}{p_{i+1}^*}}}^{[\ell,m_{i+1}]}\|\nabla \delta c\|_{\dot{\B}_{p_{i+1},1}^{\frac{d}{p_{i}}-1}}^{m_{i+1}}
\lesssim\Big(\| v_1\|_{\dot{\B}_{p,1}^{\frac dp }}^{\ell}+\sum^{R}_{\iota=i+1}\| v_1\|_{\dot{\B}_{p_\iota,1}^{\frac {d}{p_\iota}}}^{m_\iota}\Big)\| \delta c\|_{\dot{\B}_{p_{i+1},1}^{\frac{d}{p_{i+1}}}}^{m_{i+1}}
\end{eqnarray}
which gives
\begin{eqnarray}
\int^t_0\|\nabla v_1\|_{\dot{\B}_{p,1}^{\frac dp-\frac {d}{p_{i+1}^*}}}^{[\ell,m_{i+1}]}\|\nabla \delta c\|_{\dot{\B}_{p_{i+1},1}^{\frac{d}{p_{i}}-1}}^{m_{i+1}}\lesssim\int_0^t \tilde{X}(s)\delta X(s).
\end{eqnarray}
We can handle the other right-hand side terms similarly, except the last one for which we have
\begin{eqnarray}\label{gub}&&\nonumber\int^t_0\|\nabla\delta  c\|_{\dot{\B}_{p_{i+1},1}^{\frac{d}{p_{i+1}}-1}}^{[m_{i+1},h]}\|v_1\|_{\dot{\B}_{p_{i+1},1}^{\frac{d}{p_{i+1}}+1}}^{[m_{i+1},h]}\\
\nonumber&\lesssim&\int^t_0\Big(\sum^{R}_{\iota=i+1}\| \delta a\|_{\dot{\B}_{p_\iota,1}^{\frac {d}{p_\iota}}}^{m_\iota}+\| \delta \|_{\dot{\B}_{2,1}^{\frac d2 }}^{h}\Big)\Big(\sum^{R}_{\iota=i+1}\| v_1\|_{\dot{\B}_{p_\iota,1}^{\frac {d}{p_\iota}}}^{m_\iota}+\| v_1\|_{\dot{\B}_{2,1}^{\frac d2 +1}}^{h}\Big)\\
&\lesssim&\int_0^t \tilde{X}(s)\delta X(s).
\end{eqnarray}
Consequently, we reach
$$\int_0^t\sum_{j\in J_i}2^{\frac{d}{p_i}j}\|\ddj(v_1\cdot\nabla\delta c)\|_{L^{p_i}}\lesssim\int_0^t \tilde{X}(s)\delta X(s).$$
On the other hand, for the term $\delta v\cdot\nabla c_2$, we use the decomposition \eqref{S1}, \eqref{S2}-\eqref{S3}. For the medium-frequency regime $m^\iota, \iota\in[1,R]$, we have
\begin{eqnarray*}
\sum_{j\in J_i}2^{\frac{d}{p_i}j}\|\dot S_{j-1}\delta v^{m_\iota}\cdot\nabla\ddj c_{2}\|_{L^{p_i}}\leq C\|\dot S_{j-1}\delta v^{m_\iota}\|_{L^\infty}\|c_2\|_{\dot{\B}_{2,1}^{\frac{d}{2}+1}}^{h}
\lesssim\tilde{X}(s)\delta X(s)
\end{eqnarray*}
and
\begin{multline}
\sum_{j\in J_i}2^{\frac{d}{p_i}j}\|\delta\bar R_j^1(v^{m_\iota},c_2)\|_{L^{p_i}}\leq C\Bigl(\|\nabla \delta v^{m_\iota}\|_{\dot{\B}_{p,1}^{\frac dp-\frac {d}{p_{i+1}^*}}}^{[\ell,m_{i+1}]}\|\nabla  c_2\|_{\dot{\B}_{p_{i+1},1}^{\frac{d}{p_{i}}-1}}^{m_{i+1}}
+{\|\nabla \delta v^{m_\iota}\|_{L^\infty}}\|\nabla c_2\|_{\dot{\B}_{p_i,1}^{\frac{d}{p_i}-1}}^{[m_{i},h]}\hfill\cr\hfill
+
\|\nabla c_2\|_{\dot{\B}_{p,1}^{\frac dp-\frac {d}{p_{i+1}^*}}}^{[\ell,m_{i+1}]}\|\delta v^{m_\iota}\|_{\dot{\B}_{p_{i+1},1}^{\frac{d}{p_i}}}^{m_{i+1}}
+\|\nabla c_2\|_{L^\infty}\|\delta v^{m_\iota}\|_{\dot{\B}_{p_i,1}^{\frac{d}{p_i}}}^{[m_{i},h]}+ \|\nabla c_2\|_{\dot{\B}_{p_{i+1},1}^{\frac{d}{p_{i+1}}-1}}^{[m_{i+1},h]}\|\delta v^{m_\iota}\|_{\dot{\B}_{p_{i+1},1}^{\frac{d}{p_{i+1}}+1}}^{[m_{i+1},h]}\Bigl)\\
\lesssim\tilde{X}(s)\delta X(s).\end{multline}
Then, integrating in time implies
$$\int_0^t\sum_{j\in J_i}2^{\frac{d}{p_i}j}\|\ddj(\delta v\cdot\nabla c_2)\|_{L^{p_i}}\lesssim\int_0^t \tilde{X}(s)\delta X(s).$$
The remaining terms in $\delta Q^{1}_{j}$ and $\delta Q^{2}_{j}$ can be treated in a similar fashion. We obtain
\begin{eqnarray*}
\|\delta c\|_{L^\infty_T(\dot{\B}_{{p_i},1}^{\frac{d}{{p_i}}})}^{m_{i} } +\|\delta v\|_{L^\infty_T(\dot{\B}_{{p_i},1}^{\frac{d}{{p_i}}})}^{m_{i} }\leq \delta X_{0}+\int_0^t \tilde{X}(s)\delta X(s).
\end{eqnarray*}
\subsubsection{Error estimates in low frequencies}
As in the medium-frequency analysis, we have
\begin{eqnarray}
\| \delta c\|^{\ell}_{L^\infty_T(\dot{\mathbb{B}}_{p,1}^{\frac{d}{p}})}+\|\delta v\|^{\ell}_{L^\infty_T(\dot{\mathbb{B}}_{p,1}^{\frac{d}{p}})}
\leq
\delta X_{0}+\int_0^t\sum_{j\leq \J-RN_0}2^{\frac{d}{p}j}\| (\delta Q^1_{j},\delta Q^2_{j})\|_{L^p}.
\end{eqnarray}
Again, we focus on the term $\delta Q^1_j$ and specifically the term $\dot\Delta_j(v_1\cdot\nabla \delta c + \delta v \cdot \nabla c_2)$. For the first term, recalling that $\dot\Delta_j(v_1\cdot\nabla \delta c)=\dot S_{j-1}v_{1}\cdot\nabla{\delta c}_{j}+\delta R_j^1$, we have
\begin{eqnarray*}
\sum_{j\leq \J-RN_0}2^{\frac{d}{p}j}\|\ddj(v_1\cdot\nabla\delta c)\|_{L^p}\leq C(\| v_{1}\|_{\dot{\B}_{p,1}^{\frac{d}{p}}}\|\nabla\delta  c\|_{\dot{\B}_{p,1}^{\frac{d}{p}}}^\ell+\| \nabla v_{1}\|_{\dot{\B}_{p,1}^{\frac{d}{p}}}\|\delta  c\|_{\dot{\B}_{p,1}^{\frac{d}{p}}}).
\end{eqnarray*}
For the second term in $\delta Q^1_j$, we have
\begin{eqnarray*}
\sum_{j\leq \J-RN_0}2^{\frac{d}{p}j}\|\ddj(\delta v\cdot\nabla c_2)\|_{L^p}\leq C\|\delta v\|_{\dot{\B}_{p,1}^{\frac{d}{p}}}\|\nabla c_2\|_{\dot{\B}_{p,1}^{\frac{d}{p}}}
\end{eqnarray*}
and thus
\begin{eqnarray*}
\int_0^t\sum_{j\leq  \J-RN_0}2^{\frac{d}{p}j}\|\ddj(v_1\cdot\nabla\delta c+\delta v\cdot\nabla c_2\big)\|_{L^p}\leq C \int_0^t \tilde{X}(s)\delta X(s).
\end{eqnarray*}
Then, performing similar computations for the other terms in $\delta Q^1_j$ and $\delta Q^2_j$, we get
\begin{eqnarray*}
\|\delta c\|_{L^\infty_T(\dot{\B}_{{p},1}^{\frac{d}{{p}}})}^{\ell} +\|\delta v\|_{L^\infty_T(\dot{\B}_{{p},1}^{\frac{d}{{p}}})}^{\ell }\leq \delta X_{0}+\int_0^t \tilde{X}(s)\delta X(s).
\end{eqnarray*} Gathering the previous estimates, we obtain (\ref{error}).

\section{Proof of Theorem \ref{thm3}}
Let us recall the following proposition concerning perturbed solutions of the porous medium equation.
\begin{prop}[\cite{CBD3}]\label{prop1}
Let $1 \leq p < \infty$ and assume that $\mathcal{N}_0 - \bar{\mathcal{N}} \in \dot{\mathbb{B}}_{p,1}^{\frac{d}{p}}$ with $\bar{\mathcal{N}} > 0$. There exists a constant $c_0 > 0$ such that if
$$\left\|\mathcal{N}_0 - \bar{\mathcal{N}}\right\|_{\dot{\mathbb{B}}_{p,1}^{\frac dp}} \leq c_0,$$
then, the following equation
\begin{eqnarray}\label{heat}
\partial_{t}\mathcal{N}-\Delta P(\mathcal{N})=0
\end{eqnarray}
supplemented with initial data $\mathcal{N}_0$ has a unique global solution $\mathcal{N}$ such that $\mathcal{N} - \bar{\mathcal{N}} \in \mathcal{C}_b(\mathbb{R}^+; \dot{\mathbb{B}}_{p,1}^{\frac dp}) \cap L^1(\mathbb{R}^+; \dot{\mathbb{B}}_{p,1}^{\frac dp+2})$ and
\begin{eqnarray}\label{Y}
Y(t)\triangleq\|\N\|_{L^\infty_T(\dot{\B}^{\frac dp}_{p,1})}+\|\N\|_{L^1_T(\dot{\B}^{\frac dp+2}_{p,1})}\lesssim c_0.
\end{eqnarray}
\end{prop}
To prove Theorem \ref{thm3}, we use the uniform bounds established in Theorem \ref{thm2}. For $r\in[1,2]$, we have
\begin{eqnarray*}
\|\W\|_{L^r_T(\dot{\B}_{{p},1}^{\frac{d}{{p}}})}&\lesssim&\big(\|\W\|^{\ell}_{L^2_T(\dot{\B}_{{p_i},1}^{\frac{d}{{p_i}}})}+\sum^{R}_{i=1}\|\W\|^{m_i}_{L^2_T(\dot{\B}_{{p_i},1}^{\frac{d}{{p_i}}})}+\varepsilon\|(\nabla c,v)\|^{h}_{L^2_T(\dot{\B}_{{2},1}^{\frac{d}{{2}}+1})}\big)^{2-\frac{2}{r}}\\
&\quad\times&\big(\|\W\|^{\ell}_{L^1_T(\dot{\B}_{{p_i},1}^{\frac{d}{{p_i}}})}+\sum^{R}_{i=1}\|\W\|^{m_i}_{L^1_T(\dot{\B}_{{p_i},1}^{\frac{d}{{p_i}}})}+\varepsilon\|(\nabla c,v)\|^{h}_{L^1_T(\dot{\B}_{{2},1}^{\frac{d}{{2}}+1})}\big)^{\frac{2}{r}-1}\\&\lesssim\varepsilon^{\frac{2}{r}-1}.
\end{eqnarray*}
Using that $\frac{\nabla P({\rho}^\varepsilon)}{{\rho}^\varepsilon}=\check{\gamma}(c
+\bar c)\nabla c$, we obtain the desired estimates for the damped mode:
\begin{eqnarray}\label{conver1}\left\|{v}^\varepsilon+\frac{\nabla P({\rho}^\varepsilon)}{{\rho}^\varepsilon}\right\|_{L^{r}_T(\dot{\mathbb{B}}_{p,1}^{\frac{d}{p}})}=\left\|\W\right\|_{L^{r}_T(\dot{\mathbb{B}}_{p,1}^{\frac{d}{p}})}\leq C\varepsilon^{\frac{2}{r}-1}.\end{eqnarray}

Define $\delta D=\rho^\varepsilon-\N$. In light of Taylor expansion, we have, for $f\in\mathcal{S}'(\R^d)$,
$$P(f)-P(f^*)=P'(f^*)(f-f^*)+\tilde P(f^*)(f-f^*).$$
Since $\bar\rho=\bar\N$, the errr $\delta D$ fulfills
$$\partial_t \delta D-\mu\Delta \delta D=\mathcal{E}+\mathcal{S}$$
with $\mu=P'(\bar\rho)$,
$$\mathcal{E}=\Delta\big((\tilde P(\rho)-\tilde P(\bar\rho))\delta D+(\tilde P(\rho)-\tilde P(\N))(\N-\bar\N)\big)\quad\text{and}\quad \mathcal{S}=-\div(\rho \W).$$
Consequently Duhamel's formula indicates that
$$\delta D=e^{\mu\Delta t}\delta D_0+\int^{t}_{0}e^{\mu\Delta (t-s)}\big(\mathcal{E}+\mathcal{S}\big).$$
Standard estimates for the heat equation give, for $\delta>0$,
\begin{align}
\|\delta D\|_{L^\infty_T(\dot{\B}^{\frac dp-\delta}_{p,1})\cap L^{\frac{2}{1+\delta}}_T(\dot{\B}^{\frac dp+1}_{p,1})}\lesssim\|\delta D_0\|_{\dot{\B}^{\frac dp-\delta}_{p,1}}+\|\mathcal{E}\|_{L^{\frac{2}{1+\delta}}_T(\dot{\B}^{\frac dp-1}_{p,1})}+\|\mathcal{S}\|_{L^{\frac{2}{1+\delta}}_T(\dot{\B}^{\frac dp-1}_{p,1})}.
\end{align}
Hence, we are left with controlling the right-hand side terms. For the first term in $\mathcal{E}$, we have
\begin{align*}
\|\Delta\big((\tilde P(\rho)-\tilde P(\bar\rho))\delta D\big)\|_{L^{{\frac{2}{1+\delta}}}_T(\dot{\B}^{\frac dp-1}_{p,1})}&\lesssim\|c\|_{L^\infty_T(L^\infty)}\|\delta D\|_{L^{\frac{2}{1+\delta}}_T(\dot{\B}^{\frac dp+1}_{p,1})}+\|\delta D\|_{L^{\frac{2}{\delta}}_T(\dot{\B}^{\frac dp}_{p,1})}\|c\|_{L^2_T(\dot{\B}^{\frac dp+1}_{p,1})}
\\&\lesssim X(t)\|\delta D\|_{L^\infty_T(\dot{\B}^{\frac dp-\delta}_{p,1})\cap L^{\frac{2}{1+\delta}}_T(\dot{\B}^{\frac dp+1}_{p,1})}.
\end{align*}
The other term in $\mathcal{E}$ can be treated in a similar manner and we obtain
$$\|\Delta\big((\tilde P(\rho)-\tilde P(N))(\N-\bar\N)\big)\|_{L^{{\frac{2}{1+\delta}}}_T(\dot{\B}^{\frac dp-1}_{p,1})}\lesssim Y(t)\|\delta D\|_{L^\infty_T(\dot{\B}^{\frac dp-\delta}_{p,1})\cap L^{\frac{2}{1+\delta}}_T(\dot{\B}^{\frac dp+1}_{p,1})}.$$
Gathering the above estimates, we arrive at
$$\|\mathcal{E}\|_{L^{{\frac{2}{1+\delta}}}_T(\dot{\B}^{\frac dp-1}_{p,1})}\lesssim \big(X(t)+Y(t)\big)\|\delta D\|_{L^\infty_T(\dot{\B}^{\frac dp-\delta}_{p,1})\cap L^{\frac{2}{1+\delta}}_T(\dot{\B}^{\frac dp+1}_{p,1})}.$$
Concerning $\mathcal{S}$, for $\delta\in(0,1]$, there holds
\begin{eqnarray*}
\|\div(\rho \W)\|_{L^{\frac{2}{1+\delta}}_T(\dot{\B}^{\frac dp-1}_{p,1})}\lesssim(1+\|c\|_{L^\infty_T(\dot{\B}^{\frac dp}_{p,1})})\|\W\|_{L^{\frac{2}{1+\delta}}_T(\dot{\B}^{\frac dp}_{p,1})}
\lesssim\varepsilon^{\delta} X(t),
\end{eqnarray*}
where we employed (\ref{conver1}) in the last inequality. Using that $X(t),Y(t)\ll1$, we conclude that
\begin{eqnarray*}
&&\|\delta D\|_{L^\infty_T(\dot{\B}^{\frac dp-\delta}_{p,1})\cap L^r_t(\dot{\B}^{\frac dp+1}_{p,1})}\leq\|\delta D_0\|_{\dot{\B}^{\frac dp-\delta}_{p,1}}+\varepsilon^{\delta}X(t)
\end{eqnarray*}
which, provided that $\|\delta D_0\|_{\dot{\B}^{\frac dp-\delta}_{p,1}}\leq\varepsilon^{\delta}$, implies
\begin{eqnarray*}\|\delta D\|_{L^\infty_T(\dot{\B}^{\frac dp-\delta}_{p,1})}\leq\varepsilon^{\delta}.
\end{eqnarray*}
This concludes the proof of Theorem \ref{thm3}.

\section*{Acknowledgments}

T. Crin-Barat is supported by the project ANR-24-CE40-3260 – Hyperbolic Equations, Approximations $\&$ Dynamics (HEAD).

\section*{Data availability statement}

Data sharing is not applicable to this article, as no datasets were generated or analyzed during the current study.

\section*{Conflict of interest statement}

The authors declare that they have no conflict of interest.

\appendix

\section{Toolbox and proof of the product and commutator estimates}\label{Appendix}
\subsection{Toolbox}
\begin{lemma}[\cite{DecayNSCP}]\label{SimpliCarre}
Let $X : [0,T]\to \mathbb{R}_+$ be a continuous function such that $X^2$ is differentiable. Assume that there exists 
 a constant $B\geq 0$ and  a measurable function $A : [0,T]\to \mathbb{R}_+$ 
such that 
 $$\frac{1}{2}\frac{d}{dt}X^2+BX^2\leq AX\quad\hbox{a.e.  on }\ [0,T].$$ 
 Then, for all $t\in[0,T],$ we have
$$X(t)+B\int_0^tX(s)\leq X_0+\int_0^tA(s).$$
\end{lemma}

\subsection{Proof of Lemma \ref{product}}
 Employing the Bony paraproduct decomposition, we have
\begin{equation*}
\ddj( fg)=T_{f}g+R(f,g)+T_g f,
\end{equation*}
where
$$T_{f}g:=\ddj(\sum_{j'\in\Z}\dot{S}_{j'-1}f\dot{\Delta}_{j'}g),$$
$$R(f,g):=\ddj(\sum_{j'\in\Z}\tilde{\dot{\Delta}}_{j'}f\cdot \dot{\Delta}_{j'} g)$$
and
$$T_g f:=\ddj(\sum_{j'\in\Z}\dot{S}_{j'-1}g \dot{\Delta}_{j'} f),$$
where $\tilde{\dot{\Delta}}_{j'}f:=\sum\limits_{|k-j'|\leq1}{\dot{\Delta}}_{k}f$.

\underline{Case (1): High frequencies.} Let $j\geq J^\varepsilon$ and $f\in \mathcal{S}'_h(\R^d)$. For the remainder term $R(f,g)$, using the Fourier localization properties of the $\ddj$ dyadic block, there exists a positive $N_1$ such that
\begin{multline}\label{support}
\mathcal{F}(\sum_{j'\in\Z}\ddj(\tilde{\dot{\Delta}}_{j'}f\dot{\Delta}_{j'}g))(\xi)=\sum_{j'\in\Z}\int_{\R^d}\varphi(\frac{|\xi|}{2^j})
\tilde{\varphi}(\frac{|\xi-\eta|}{2^{j'}})\varphi(\frac{|\eta|}{2^{j'}})\hat f(\xi-\eta)\widehat{g}(\eta)d\eta\\
=\sum_{j'\geq j-N_1}\int_{\R^d}\varphi(\frac{|\xi|}{2^j})
\tilde{\varphi}(\frac{|\xi-\eta|}{2^{j'}})\varphi(\frac{|\eta|}{2^{j'}})\hat f(\xi-\eta)\widehat{g}(\eta)d\eta.
\end{multline}
Indeed, due to the definition of $\varphi$ and $\tilde\varphi$, we find that
$$|\xi|\in[\frac{3}{4}2^j,\frac{8}{3}2^j],\quad
|\xi-\eta|\in[\frac{3}{8}2^{j'},\frac{16}{3}2^{j'}]\quad\text{and}\quad
|\eta|\in[\frac{3}{4}2^{j'},\frac{8}{3}2^{j'}].$$
Now for a fixed $j$, if $j'<j-N_1$, then
$$|\xi|\leq\big||\xi-\eta|+|\eta|\big|\leq2^{j'+2N_1}.$$
According to the support of  $|\xi|$, provided we select $N_1$ large enough, we obtain (\ref{support}).
Hence, we have
\begin{align*}
\sum_{j\geq J}2^{sj}\|\ddj(\sum_{j'\in\Z}\tilde{\dot{\Delta}}_{j'}f \dot{\Delta}_{j'} g)\|_{L^2}&\lesssim
\sum_{j\geq J}2^{(s+\frac{2d}{p_1}-\frac{d}{2})j}\|\ddj(\sum_{j'\in\Z}\tilde{\dot{\Delta}}_{j'}f \dot{\Delta}_{j'} g)\|_{L^\frac{p_1}{2}}\\
&\lesssim\sum_{j\geq J}2^{(s+\frac{2d}{p_1}-\frac{d}{2})j}\sum_{j'\geq j-N_1}\|\tilde{\dot{\Delta}}_{j'}f\|_{L^{p_1}}\|\dot{\Delta}_{j'} g\|_{L^{p_1}}\\
&=\sum_{j\geq J}\left(\sum_{j'\geq j-N_1}2^{(s+\frac{2d}{p_1}-\frac{d}{2})(j-j')}2^{(s+\frac{2d}{p_1}-\frac{d}{2})j'}\|\tilde{\dot{\Delta}}_{j'}f\|_{L^{p_1}}\|\dot{\Delta}_{j'} g\|_{L^{p_1}}\right).
\end{align*}
If $s+\frac{2d}{p_1}-d>0$, we can select a $N_0$ satisfying  $N_0\geq N_1+1$ and apply Young's inequality to get, for $s-\frac{d}{p}=s_1+s_2-\frac{2d}{p_1}$,
\begin{align*}
\sum_{j\geq J}2^{sj}\|\ddj(\sum_{j'\in\Z}\tilde{\dot{\Delta}}_{j'}f\dot{\Delta}_{j'}g)\|_{L^2}
&\lesssim\sum_{j\geq J}\left(\sum_{j'\geq j-N_0-1}2^{(s+\frac{2d}{p_1}-\frac{d}{2})(j-j')}2^{(s+\frac{2d}{p_1}-\frac{d}{2})j'}\|\tilde{\dot{\Delta}}_{j'}f\|_{L^{p_1}}\|\dot{\Delta}_{j'} g\|_{L^{p_1}}\right)\\
&\lesssim
\sum_{j\geq J-N_0-1}2^{s_{1}j}\|{\tilde{\dot{\Delta}}}_{j}f\|_{L^{p_1}}\sup_{j\geq J-N_0-1}2^{s_{2}j}\|\dot{\Delta}_{j} g\|_{L^{p_1}}\\
&\lesssim
\sum_{j\geq J-N_0}2^{s_{1}j}\| {\dot{\Delta}}_{j}f\|_{L^{p_1}}\sup_{j\geq J-N_0}2^{s_{2}j}\|\dot{\Delta}_{j} g\|_{L^{p_1}}\\&
\lesssim\|f\|_{\dot{\B}_{p_1,1}^{s_{1}}}^{[m_1,h]}\|g\|_{\dot{\B}_{p_1,1}^{s_{2}}}^{[m_2,h]}.
\end{align*}
Concerning the paraproduct term $T_fg$, there exists a $N_2\in \mathbb{Z}$ such that
\begin{multline}\label{support2}
\mathcal{F}(\sum_{j'\in\Z}\ddj(\dot{S}_{j'-1}f\dot{\Delta}_{j'}g))(\xi)=\sum_{j'\in\Z}\int_{\R^d}\varphi(\frac{|\xi|}{2^j})
{\chi}(\frac{|\xi-\eta|}{2^{j'-1}})\varphi(\frac{|\eta|}{2^{j'}})\hat f(\xi-\eta)\widehat{g}(\eta)d\eta\\
=\sum_{|j-j'|\leq N_2}\int_{\R^d}\varphi(\frac{|\xi|}{2^j})
{\chi}(\frac{|\xi-\eta|}{2^{j'-1}})\varphi(\frac{|\eta|}{2^{j'}})\hat f(\xi-\eta)\widehat{g}(\eta)d\eta.
\end{multline}
The proof of the above equality is similar. to that of (\ref{support}). The triangle inequality ensures that
$$|\xi|\in[\frac{3}{4}2^j,\frac{8}{3}2^j],\,\,
|\xi-\eta|\in[0,\frac{4}{3}2^{j'}]\quad \text{and}\quad
|\eta|\in[\frac{3}{4}2^{j'},\frac{8}{3}2^{j'}]\Rightarrow|\xi|\in[\frac{4}{3}2^{j'},2^{j'+2}].$$
Then, choosing $N_2$ sufficiently large leads to a contradiction, ensuring (\ref{support2}). Naturally  (\ref{support2}) implies
$$\sum_{j\geq J}\ddj(\sum_{j'\in\Z}\dot{S}_{j'-1}f\dot{\Delta}_{j'}g)\sim
\sum_{j\geq J-N_2}\dot{S}_{j-1}f\dot{\Delta}_{j}g$$
and
\begin{multline*}
\|\sum_{j\geq J}2^{sj}\ddj(\sum_{j'\in\Z}\dot{S}_{j'-1}f\dot{\Delta}_{j'}g)\|_{L^2}\\
\lesssim\sum_{j\in [J-N_2,J)}2^{sj}\|\dot{S}_{j-1}f\|_{L^{p^{*}_{1}}}\|\dot{\Delta}_{j}g\|_{L^{p_1}}
+\sum_{j\geq J}2^{sj}\|\dot{S}_{j-1}f\|_{L^\infty}\|\dot{\Delta}_{j}g\|_{L^2}.
\end{multline*}
Employing Sobolev embeddings, we obtain
$$\sup_{j\in [J-N_2,J)}\|\dot{S}_{j-1}f\|_{L^{p^{*}_{1}}}\lesssim \|f\|_{\dot{\B}_{p,1}^{\frac dp-\frac {d}{p_1^*}}}^{[\ell,m_1]},$$
$$\sup_{j\geq J}\|\dot{S}_{j-1}f\|_{L^\infty}\lesssim \|f\|_{L^\infty}\quad\text{and}\quad\sum_{j\geq J}2^{sj}\|\dot{\Delta}_{j} g\|_{L^{2}}\lesssim\|g\|^{h}_{\dot{\B}^{s}_{2,1}},$$
provided that $p_1\geq\frac{2p}{p-2}$. Moreover, by further selecting $N_0\geq N_2$, it holds that
$$\sum_{j\in [J-N_2,J)}2^{sj}\|\dot{\Delta}_{j}g\|_{L^{p_1}}\lesssim\sum_{j\in [J-N_0,J)}2^{sj}\|\dot{\Delta}_{j}g\|_{L^{p_1}}\lesssim\|g\|^{m_1}_{\dot{\B}^{s}_{p_1,1}}.$$
The other paraproduct term can be handled symmetrically. This concludes the high-frequency estimates.

\underline{Case (2): Products in medium frequencies.} Recall that $J_i=[J-N_0 i,J-(N_0-1)i[$. 
For the remainder term, proceeding as in the high-frequency regime, we find that there exists a $N_1$ such that
\begin{align*}
\sum_{j\in J_i}2^{sj}\|\ddj(\sum_{j'\in\Z}\tilde{\dot{\Delta}}_{j'}f \dot{\Delta}_{j'} g)\|_{L^{p_i}}\lesssim
\sum_{j\in J_i}2^{(s+\frac{2d}{p_{i+1}}-\frac{d}{p_{i}})j}\|\ddj(\sum_{j'\in\Z}\tilde{\dot{\Delta}}_{j'}f \dot{\Delta}_{j'}g)\|_{L^\frac{p_{i+1}}{2}}\\
\lesssim\sum_{j\in J_i}2^{(s+\frac{2d}{p_{i+1}}-\frac{d}{p_{i}})j}\sum_{j'\geq j-N_1}\|\tilde{\dot{\Delta}}_{j'}f\|_{L^{p_{i+1}}}\|\dot{\Delta}_{j'}g\|_{L^{p_{i+1}}}\\
=\sum_{j\in J_i}\left(\sum_{j'\geq j-N_1}2^{(s+\frac{2d}{p_{i+1}}-\frac{d}{p_{i}})(j-j')}2^{(s+\frac{2d}{p_{i+1}}-\frac{d}{p_{i}})j'}\|\tilde{\dot{\Delta}}_{j'}f\|_{L^{p_{i+1}}}\|\dot{\Delta}_{j'}g\|_{L^{p_{i+1}}}\right).
\end{align*}
Since $J_i=[J-N_0 i,J-(N_0-1)i[$, we can select a $N_0$ satisfying $N_0\geq N_1+1$, then, if $s+\frac{2d}{p_{i+1}}-\frac{d}{p_{i}}>0$, Young's inequality gives
\begin{multline*}
\sum_{j\in J_i}2^{sj}\|\ddj(\sum_{j'\in\Z}{\dot{\Delta}}_{j'}f\dot{\Delta}_{j'}g)\|_{L^{p_i}}\lesssim
\sum_{j\geq J-N_0(i+1) i}2^{s_1j}\|\tilde{\dot{\Delta}}_{j}f\|_{L^{p_{i+1}}}\\ \times\sup_{j\geq J-N_0 (i+1)}2^{s_2j}\|\dot{\Delta}_{j} g\|_{L^{p_{i+1}}}
\lesssim  \|f\|_{\dot{\B}_{p_{i+1},1}^{s_1}}^{[m_{i+1},h]}  \|g\|_{\dot{\B}_{p_{i+1},1}^{s_2}}^{[m_{i+1},h]}
\end{multline*}
with $s_1+s_2=s+\frac{2d}{p_{i+1}}-\frac{d}{p_{i}}$ and this concludes the analysis of the remainder term.

For the paraproduct terms, as in the high-frequency regime, we can choose a $N_2$ such that
\begin{multline*}
\|\sum_{j\in J_i}2^{sj}\ddj(\sum_{j'\in\Z}\dot{S}_{j'-1}f\dot{\Delta}_{j'}g)\|_{L^{p_i}}
\lesssim
\sum_{j\geq J-N_0 i-N_2}2^{sj}\|\dot{S}_{j-1}f\dot{\Delta}_{j}g\|_{L^{p_i}}\\
\lesssim\sum_{j\in [J-N_0 i-N_2,J-N_0 i)}2^{sj}\|\dot{S}_{j-1}f\|_{L^{p^{*}_{i+1}}}\|\dot{\Delta}_{j}g\|_{L^{p_{i+1}}}
+\sum_{j\geq J-N_0 i}2^{sj}\|\dot{S}_{j-1}f\|_{L^\infty}\|\dot{\Delta}_{j}g\|_{L^{p_i}}.
\end{multline*}
Since
$$\sup_{j\in [J-N_0i-N_2,J-N_0 i)}\|\dot{S}_{j-1}f\|_{L^{p^{*}_{i+1}}}\lesssim \|f\|_{\dot{\B}_{p,1}^{\frac dp-\frac {d}{p_{i+1}^*}}}^{[\ell,m_{i+1}]},$$
$$\sup_{j\geq J-N_0 i}\|\dot{S}_{j-1}f\|_{L^\infty}\lesssim \|f\|_{L^\infty},\quad\sum_{j\geq J-N_0 i}2^{sj}\|\dot{\Delta}_{j} g\|_{L^{2}}\lesssim\|g\|^{[m_i,h]}_{\dot{\B}^{s}_{p_i,2,1}},$$
if we choose $N_0\geq N_2$, there holds
$$\sum_{j\in [J-N_0 i-N_2,J-N_0 i)}2^{sj}\|\dot{\Delta}_{j} g\|_{L^{p_{i+1}}}\lesssim\sum_{j\in [J-N_0(i+1),J-N_0 i)}2^{sj}\|\dot{\Delta}_{j} g\|_{L^{p_{i+1}}}\lesssim\|g\|^{m_i}_{\dot{\B}^{s}_{p_{i+1},1}},$$
which gives the desired estimate in the medium-frequency regime and concludes the proof of Lemma \ref{product}.\qed

\subsection{Proof of Lemma \ref{commutator}}
Recall that $\mathfrak{R}_j:= \dot S_{j-1}f\,\ddj g-\ddj(fg)$. By Bony decomposition, we have
\begin{equation*}
\mathfrak{R}_j:=R_{1,j}+R_{2,j}+R_{3,j}
\end{equation*}
where
$$R_{1,j}=\ddj(\sum_{j'\in\Z}\tilde{\dot{\Delta}}_{j'}w\cdot \dot{\Delta}_{j'} z), \quad R_{2,j}=\ddj(\sum_{j'\in\Z}\dot{\Delta}_{j'}w\cdot\dot{S}_{j'-1}z)$$
and
$$R_{3,j}=\ddj(\sum_{j'\in\Z}\dot{S}_{j'-1}w\cdot \dot{\Delta}_{j'} z)-\dot{S}_{j-1}w\cdot\ddj z.$$

\underline{Case (1): High frequencies.} Let $j\geq J$. 
The terms $R_{1,j}$ and $R_{2,j}$ can be handled as in the proof of Lemma \ref{product}, so we focus on $R_{3,j}$. Using the Fourier localization properties of the dyadic blocks, there exists an integer $N_2>0$ such that
\begin{multline}\label{Comm}
\ddj(\sum_{j'\in\Z}\dot{S}_{j'-1}w\dot{\Delta}_{j'}z)=
\sum_{|j-j'|\leq N_2}[\ddj,\dot{S}_{j'-1}w]\dot{\Delta}_{j'}z+\sum_{|j-j'|\leq 1}\left(\dot{S}_{j'-1} w-\dot{S}_{j-1} w\right)\ddj \dot{\Delta}_{j'}z.
\end{multline}
For $j'=J-1$ and choosing a $N_0$ sufficiently large, we have that the above quantity is non-zero only for $j\in[J,J+N_0]$ since $|j-j'|\leq N_2$.
Applying \cite[Lemma 2.97]{HJR} leads to
\begin{multline*}
\|\sum_{j\geq J}2^{sj}[\ddj,\dot{S}_{J-2}  w]\dot{\Delta}_{J-1}z\|_{L^2}
\lesssim\sum_{j\in[J,J+N_0]}2^{sj}\sum_{k\leq J-3}\|[\ddj,\dot{\Delta}_{k}  w]\dot{\Delta}_{J-1} z\|_{L^2}\\
\lesssim\sum_{k\leq J-3}\|\dot{\Delta}_{k}\nabla w\|_{L^{p^{*}_{1}}}\| \dot{\Delta}_{J-1} z\|_{L^{p_1}}
\lesssim\|\nabla w\|_{\dot{\B}_{p^{*}_{1},1}^{0}}^{[\ell,m_1]}\|z\|^{m_1}_{\dot{\B}^{s}_{p_1,1}}.
\end{multline*}
Similarly, one can deal with the case $j'\leq J-2$. For $j'\geq J$, there holds
\begin{multline*}
\|\sum_{j\geq J}2^{sj}\sum_{|j-j'|\leq N_2,j'\geq J}[\ddj,\dot{S}_{j'-1}  w]\dot{\Delta}_{j'}z\|_{L^2}
\lesssim\sum_{j\geq J}2^{(s-1)j}\sum_{|j-j'|\leq N_2,j'\geq J}\|\dot{S}_{j'-1}\nabla w\|_{L^\infty}\| \dot{\Delta}_{j'} z\|_{L^2}\\
\lesssim\|\nabla w\|_{\dot{\B}^{0}_{\infty,1}}\|z\|^{h}_{\dot{\B}^{s-1}_{2,1}}.
\end{multline*}
Finally, we have
\begin{multline*}
\|\sum_{j\geq J}2^{sj}\sum_{|j-j'|\leq1}\left(\dot{S}_{j'-1} w-\dot{S}_{j-1} w\right)\ddj \dot{\Delta}_{j'}z\|_{L^2}
\lesssim\sum_{j\geq J}2^{sj}\|\dot{S}_{j'-1}  w-\dot{S}_{j-1}  w\|_{L^\infty}\|\ddj  z\|_{L^2}\\
\lesssim\sup_{j\in\Z}2^j\|\ddj w\|_{L^\infty}\|z\|^{h}_{\dot{\B}^{s-1}_{2,1}}\lesssim\|\nabla w\|_{\dot{\B}^{0}_{\infty,1}}\|z\|^{h}_{\dot{\B}^{s-1}_{2,1}},
\end{multline*}
which concludes the high-frequency part.
\medbreak
\underline{Case (2): Commutators in medium frequencies.} Let $j\in J_i$. Again, we only focus $R_{3,j}$. Recalling \eqref{Comm}, for $j'=J-N_0i-1$, applying of \cite[Lemma 2.97]{HJR} leads to
\begin{multline*}
\|\sum_{j\in J_i}2^{sj}[\ddj,\dot{S}_{J-4i-2}  w]\dot{\Delta}_{J-N_0 i-1}z\|_{L^{p_i}}
\lesssim\sum_{j\in J_i}2^{sj}\sum_{k\leq J-N_0 i-3}\|[\ddj,\dot{\Delta}_{k}  w]\dot{\Delta}_{J-N_0 i-1} z\|_{L^{p_i}}\\
\lesssim\sum_{k\leq J-N_0 i-3}\|\dot{\Delta}_{k}\nabla w\|_{L^{p^*_{i+1}}}\| \dot{\Delta}_{J-N_0 i-1} z\|_{L^{p_{i+1}}}
\lesssim\|\nabla w\|_{\dot{\B}_{p^*_{i+1},1}^{0}}^{[\ell,m_{i+1}]}\|z\|^{m_{i+1}}_{\dot{\B}^{s}_{p_{i+1},1}}.
\end{multline*}
Similarly, one can treat the case $j'\leq J-N_0i-2$. For $j'\geq J-N_0 i$, there holds
\begin{multline*}
\|\sum_{j\in J_i}2^{sj}\sum_{|j-j'|\leq N_2,j'\geq J-N_0 i}[\ddj,\dot{S}_{j'-1}  w]\dot{\Delta}_{j'}z\|_{L^{p_i}}\\
\lesssim\sum_{j\in J_i}2^{(s-1)j}\sum_{|j-j'|\leq N_2,j'\geq J-N_0 i}\|\dot{S}_{j'-1}\nabla w\|_{L^\infty}\| \dot{\Delta}_{j'} z\|_{L^{p_i}}
\lesssim\|\nabla w\|_{\dot{\B}^{0}_{\infty,1}}\|z\|^{[m_i,h]}_{\dot{\B}^{s-1}_{p_i,1}}.
\end{multline*}
Finally, we have
\begin{multline*}
\|\sum_{j\in J_i}2^{sj}\sum_{|j-j'|\leq1}\left(\dot{S}_{j'-1} w-\dot{S}_{j-1} w\right)\ddj \dot{\Delta}_{j'}z\|_{L^{p_i}}
\lesssim\sum_{j\in J_i}2^{sj}\|\dot{S}_{j'-1}  w-\dot{S}_{j-1}  w\|_{L^\infty}\|\ddj  z\|_{L^{p_i}}\\
\lesssim\sup_{j\in\Z}2^j\|\ddj w\|_{L^\infty}\|z\|^{[m_i,h]}_{\dot{\B}^{s-1}_{p_i,1}}\lesssim\|\nabla w\|_{\dot{\B}^{0}_{\infty,1}}\|z\|^{[m_i,h]}_{\dot{\B}^{s-1}_{p_i,1}}.
\end{multline*}
which concludes the proof of the Lemma \ref{commutator}. \qed

\bibliographystyle{plain}
\bibliography{main.bib}

\vspace{2cm}
Timothée Crin-Barat \hfill\break\indent
{\sc Université de Toulouse, Institut de Mathématiques de Toulouse, Route de Narbonne 118, 31062 CEDEX 9 Toulouse, France, \hfill\break\indent
{\it Email address}: {\tt timothee.crin-barat@math.univ-toulouse.fr}}

\bigbreak
Zihao Song \hfill\break\indent
{\sc Mathematics and Key Laboratory of Mathematical MIIT, Nanjing University of Aeronautics and
Astronautics, Nanjing, 211106, P. R. China, \hfill\break\indent
{\it Email address}: {\tt szh1995@nuaa.edu.cn; songzh19950504@gmail.com}}

\end{document}